\documentclass{amsart}
\usepackage{amscd,amssymb,amsxtra,diagrams}
\usepackage{pstcol}

\theoremstyle{plain}
\newtheorem{Thm}{Theorem}[section]
\newtheorem{Prop}[Thm]{Proposition}
\newtheorem{Lem}[Thm]{Lemma}
\newtheorem{Cor}[Thm]{Corollary}
\newtheorem{Claim}[Thm]{Claim}

\theoremstyle{definition}
\newtheorem{Def}[Thm]{Definition}

\theoremstyle{remark}

\numberwithin{equation}{section}

\newcommand{\res}{\upharpoonright}
\newcommand{\cf}{\operatorname{cf}}
\newcommand{\ci}{\operatorname{coi}}
\newcommand{\ran}{\operatorname{ran}}
\newcommand{\Con}{\operatorname{Con}}
\newcommand{\dom}{\operatorname{dom}}

\newcommand{\st}{\operatorname{st}}

\newcommand{\inv}{{-1}}
\newcommand{\RR}{\mathbb R}
\newcommand{\QQ}{\mathbb Q}

\newcommand{\ZZ}{\mathbb Z}
\newcommand{\QR}{\mathbb Q_\mathbb R}
\newcommand{\cI}{\mathcal I}
\newcommand{\cF}{\mathcal F}
\newcommand{\cD}{\mathcal D}
\newcommand{\cR}{\mathcal R}
\newcommand{\cA}{\mathcal A}
\newcommand{\cO}{\mathcal O}

\newcommand{\uG}{\mathbf G}
\newcommand{\uS}{\mathbf S}
\newcommand{\uN}{\mathbf N}

\newcommand{\eps}{\varepsilon}
\newcommand{\aff}{{\mathrm{aff}}}
\newcommand{\cvx}{\mathrm{cvx}}

\begin{document}


\title{Asymptotic cones of finitely presented groups}


\author{Linus Kramer}
\address
{Fachbereich Mathematik\\
TU Darmstadt \\
Schlossgartenstr. 7 \\
D-64289 Darmstadt \\
Germany}
\email{kramer@mathematik.tu-darmstadt.de}
\thanks{The research of the first and third authors
was supported by Heisenberg Fellowships of the DFG.
The research of the second author was supported
by the Israel Science Foundation. This paper is
number 818 in the cumulative list of the second author's
publications. The research of the fourth author was partially 
supported by NSF Grants.} 

\author{Saharon Shelah}
\address
{Mathematics Department \\
The Hebrew University \\
Jerusalem \\
Israel}
\email{shelah@math.huji.ac.il}

\author{Katrin Tent}
\address
{Mathematisches Institut \\
Universit\"{a}t W\"{u}rzburg \\
Am Hubland \\
D-97074 W\"{u}rzburg \\
Germany}
\email{tent@mathematik.uni-wuerzburg.de}

\author{Simon Thomas}
\address
{Mathematics Department \\
Rutgers University \\
110 Frelinghuysen Road \\
Piscataway \\
New Jersey 08854-8019 \\
USA}
\email{sthomas@math.rutgers.edu}


\begin{abstract}
Let $G$ be a connected semisimple Lie group with at least
one absolutely simple factor $S$ such that 
$\mathbb{R}\text{-rank}(S) \geq 2$ and let $\Gamma$ be a
uniform lattice in $G$.
\begin{enumerate}
\item[(a)] If $CH$ holds, then $\Gamma$ has a unique
asymptotic cone up to homeomorphism.
\item[(b)] If $CH$ fails, then $\Gamma$ has 
$2^{2^{\omega}}$ asymptotic cones up to homeomorphism.
\end{enumerate}
\end{abstract}

\maketitle

\section{Introduction} \label{S:intro}

Let $\Gamma$ be a finitely generated group equipped with a fixed finite
generating set and let $d$ be the corresponding word metric. 
Consider the sequence of metric spaces
$X_{n} = (\Gamma, d_{n})$ for $n \geq 1$, where 
$d_{n}(g,h) = d(g,h)/n$. In \cite{g1}, Gromov proved that if
$\Gamma$ has polynomial growth, then the sequence 
$( X_{n} \mid n \geq 1 )$ of metric spaces converges in the pointed 
Gromov-Hausdorff topology to a complete geodesic space
$\Con_{\infty}(\Gamma)$, the asymptotic cone of $\Gamma$.
In \cite{dw}, van den Dries and Wilkie generalised the construction
of asymptotic cones to arbitrary finitely generated groups.
However, their construction involved the choice of a
nonprincipal ultrafilter $\mathcal{D}$ over the set $\omega$
of natural numbers, and it was initially not clear whether the resulting
asymptotic cone $\Con_{\mathcal{D}}(\Gamma)$ depended on the 
choice of the ultrafilter $\mathcal{D}$. In \cite{tv}, 
answering a question of Gromov \cite{g2}, Thomas and Velickovic 
constructed an example of a finitely generated group $\Gamma$ and 
two nonprincipal ultrafilters $\mathcal{A}$, $\mathcal{B}$ such that 
the asymptotic cones $\Con_{\mathcal{A}}(\Gamma)$ and 
$\Con_{\mathcal{B}}(\Gamma)$ 
were not homeomorphic. But this still left open the 
interesting question of whether there exists a finitely presented
group with more than one asymptotic cone up to 
homeomorphism. It seems almost certain that such a group 
exists; and, in fact, it seems natural to conjecture that there
exists a finitely presented group with
$2^{2^{\omega}}$ asymptotic cones up to homeomorphism.
(Recall that there are exactly $2^{2^{\omega}}$ distinct 
nonprincipal ultrafilters over $\omega$.) The main result
of this paper provides a confirmation of this 
conjecture, under the assumption that $CH$ fails.

Suppose that $G$ is a connected semisimple Lie group and let
$\Gamma$ be a uniform lattice in $G$; i.e. a
discrete subgroup such that $G/\Gamma$ is compact. (For the
existence of such a subgroup $\Gamma$, see Borel \cite{b2}.)
Then it is well-known that $\Gamma$ is finitely presented.
(For example, see \cite[Chapter 3]{vgs}, \cite[Chapter V]{h}.)
Furthermore, $\Gamma$ is quasi-isometric to $G$; and hence for
each ultrafilter $\mathcal{D}$, the asymptotic cones
$\Con_{\mathcal{D}}(\Gamma)$ and $\Con_{\mathcal{D}}(G)$ are
homeomorphic. In \cite[Section $2.2.B_{1}$]{g2}, Gromov 
suggested that `` it seems that these groups $G$ have pretty
looking finite dimensional cones $\Con_{\mathcal{D}}(G)$ which
are (essentially) independent of the choice of the 
ultrafilter $\mathcal{D}$.'' It turns out that the situation
is more interesting.

\begin{Thm} \label{T:main}
Suppose that $G$ is a connected semisimple Lie group with at least
one absolutely simple factor $S$ such that 
$\mathbb{R}\text{-rank}(S) \geq 2$ and let $\Gamma$ be a
uniform lattice in $G$.
\begin{enumerate}
\item[(a)] If $CH$ holds, then $\Gamma$ has a unique
asymptotic cone up to homeomorphism.
\item[(b)] If $CH$ fails, then $\Gamma$ has 
$2^{2^{\omega}}$ asymptotic cones up to homeomorphism.
\end{enumerate}
\end{Thm}

Here $CH$ denotes the {\em Continuum Hypothesis\/}; i.e.
the statement that $2^{\omega} = \omega_{1}$. Of course,
it is well-known that $CH$ can neither be proved nor
disproved using the usual $ZFC$ axioms of set theory. (For
example, see Jech \cite{j}.) 

In the remainder of this section, we shall sketch the main points 
of the proof for the special case when $G = SL_{m}(\mathbb{R})$ for 
some $m \geq 3$. In particular, we shall explain the unexpected
appearance of the Robinson field $\,^{\rho}\mathbb{R}_{\mathcal{D}}$
as a topological invariant of the asymptotic cone
$\Con_{\mathcal{D}}(\Gamma)$. (As we shall explain below, the
Robinson field $\,^{\rho}\mathbb{R}_{\mathcal{D}}$ is a valued
field closely related to the corresponding field 
$\,^{*}\mathbb{R}_{\mathcal{D}}$ of nonstandard real numbers.)
We shall begin by recalling the definition of an asymptotic cone of 
an arbitrary metric space.

\begin{Def} \label{D:ultra}
A {\em nonprincipal ultrafilter\/} over the set $\omega$ of
natural numbers is a collection
$\mathcal{D}$ of subsets of $\omega$ satisfying the
following conditions:
\begin{enumerate}
\item[(i)] If $A$, $B \in \mathcal{D}$, then 
$A \cap B \in \mathcal{D}$.
\item[(ii)] If $A \in \mathcal{D}$ and 
$A \subseteq B \subseteq \omega$, then
$B \in \mathcal{D}$.
\item[(iii)] For all $A \subseteq \omega$, either
$A \in \mathcal{D}$ or 
$\omega \smallsetminus A \in \mathcal{D}$.
\item[(iv)] If $F$ is a finite subset of $\omega$,
then $F \notin \mathcal{D}$.
\end{enumerate}
\end{Def}

Equivalently, if $\mu : \mathcal{P}(\omega) \to \{0,1\}$ is the function 
such that $\mu (A) = 1$ if and only if $A \in \mathcal{D}$, then $\mu$
is a finitely additive probability measure on 
$\omega$ such that $\mu (F) = 0$ for all finite
subsets $F$ of $\omega$. It is easily checked that if
$\mathcal{D}$ is a nonprincipal ultrafilter and
$( r_{n} )$ is a bounded sequence of real numbers, then there
exists a unique real number $\ell$ such that
\[
\{ n \in \omega \mid | r_{n} - \ell | < \varepsilon \}
\in \mathcal{D}
\]
for all $\varepsilon > 0$. We write
$\ell = \lim_{\mathcal{D}}r_{n}$.

\begin{Def} \label{D:cone}
Suppose that $\mathcal{D}$ is a nonprincipal ultrafilter over $\omega$.
Let $\left( X,d \right)$ be a metric space and for each $n \geq 1$,
let $d_{n}$ be the rescaled metric defined by
$d_{n}(x,y) = d(x,y)/n$. Let $e \in X$ be a fixed base point. 
Then $X_{\infty}$ is the set of all sequences $(x_{n})$ of 
elements of $X$ such that there exists a constant $c$ with
\[
d_{n}(x_{n}, e) \leq c
\]
for all $n \geq 1$. Define an equivalence relation $\sim$ on
$X_{\infty}$ by
\[
(x_{n}) \sim (y_{n}) \quad \text{ if and only if } \quad
\lim\nolimits_{\mathcal{D}}d_{n}(x_{n},y_{n}) =0;
\]
and for each $(x_{n}) \in X_{\infty}$, let $(x_{n})_{\mathcal{D}}$
be the corresponding equivalence class.

Then the asymptotic cone of $X$ is
\[
\Con_{\mathcal{D}}(X) = \{ (x_{n})_{\mathcal{D}} \mid
(x_{n}) \in X_{\infty} \}
\]
endowed with the metric
\[
d_{\mathcal{D}}((x_{n})_{\mathcal{D}}, (y_{n})_{\mathcal{D}})
= \lim\nolimits_{\mathcal{D}}d_{n}(x_{n},y_{n}).
\]
\end{Def}

If $\left( X,d \right)$ and $\left( X^{\prime},d^{\prime} \right)$ 
are metric spaces, then a map $f: X \to X^{\prime}$ is a 
{\em quasi-isometry\/} iff there exist constants $L \geq 1$ and 
$C \geq 0$ such that for all $x$, $y \in X$
\begin{itemize}
\item $\frac{1}{L}d(x,y) -C \leq d^{\prime}(f(x),f(y)) \leq
Ld(x,y) + C$;
\end{itemize}
and for all $z \in X^{\prime}$
\begin{itemize}
\item $d^{\prime}(z, f[X]) \leq C$.
\end{itemize}
The following result is well-known. (For example, see 
Proposition 2.4.6 of Kleiner-Leeb \cite{kl}.)

\begin{Prop}
Suppose that $\mathcal{D}$ is a nonprincipal ultrafilter over $\omega$ and
that $f: X \to X^{\prime}$ is a quasi-isometry of nonempty metric spaces.
Then $f$ induces a bilipschitz homeomorphism
$\varphi: \Con_{\mathcal{D}}(X) \to \Con_{\mathcal{D}}(X^{\prime})$
defined by
$\varphi( (x_{n})_{\mathcal{D}} ) = ( f(x_{n}) )_{\mathcal{D}}$.
\end{Prop}

From now on, fix some $m \geq 3$ and let $\Gamma$ be a uniform lattice in
$SL_{m}(\mathbb{R})$; i.e. $\Gamma$ is a discrete subgroup of
$SL_{m}(\mathbb{R})$ such that $SL_{m}(\mathbb{R})/\Gamma$ is compact.
(For example, let $K = k( \sqrt{\varepsilon} )$, where 
$\varepsilon = 1 + \sqrt{2}$ and $k = \mathbb{Q}(\sqrt{2})$.
Let $\sigma$ be the automorphism of $K$ over $k$ such that
$\sigma ( \sqrt{\varepsilon} ) = - \sqrt{\varepsilon}$ and let
$f$ be the Hermitian form defined by
\[
f(x,y) = x_{1}y_{1}^{\sigma} + \cdots + x_{m}y_{m}^{\sigma}. 
\]
Finally let $J$ be the ring of algebraic integers in $K$. By
\cite[Chapter 3]{vgs}, $SU(f,J)$ is a uniform lattice in
$SL_{m}(\mathbb{R})$.)
By Theorem IV.23 \cite{h}, $\Gamma$ is quasi-isometric to 
$SL_{m}(\mathbb{R})$, viewed as a metric space with respect to some
left-invariant Riemannian metric. Thus if $\mathcal{D}$ is a
nonprincipal ultrafilter over $\omega$, then $\Con_{\mathcal{D}}(\Gamma)$ 
is bilipschitz homeomorphic to $\Con_{\mathcal{D}}(SL_{m}(\mathbb{R}))$.
However, instead of working directly with the Riemannian manifold
$SL_{m}(\mathbb{R})$, it turns out to be more convenient to work
with the corresponding symmetric space, obtained by factoring out 
the maximal compact subgroup $SO_{m}(\mathbb{R})$. In more detail, recall that 
$SL_{m}(\mathbb{R})$ acts transitively as a group of isometries on 
the symmetric space $P(m, \mathbb{R})$ of positive-definite symmetric 
$m \times m$ matrices with determinant 1, via the action
\[
g \cdot A = g A g^{t}.
\]
Clearly the stabiliser of the identity matrix $I$ under this action
is the subgroup $SO_{m}(\mathbb{R})$; and since $SO_{m}(\mathbb{R})$ is compact, it
follows that $\Gamma$ also acts cocompactly on $P(m, \mathbb{R})$.
(For example, see \cite[Chapter 1]{vgs}.)
Hence, applying Theorem IV.23 \cite{h} once again, it follows that $\Gamma$ 
is also quasi-isometric to the symmetric space $P(m, \mathbb{R})$.
The invariant Riemannian metric $d$ on $P(m, \mathbb{R})$ can be described
as follows. First recall that if $A$, $B \in P(m, \mathbb{R})$, then
there exists $g \in SL_{m}(\mathbb{R})$ such that $g A g^{t}$ and
$g B g^{t}$ are simultaneously diagonal. Hence it suffices to
consider the case when 
\[
A = \Bigl(a_{i,j} \Bigr) \quad \text{ and } \quad
B = \Bigl(b_{i,j} \Bigr)
\]
are both diagonal matrices; in which case,
\[
d(A,B) = \sqrt{\sum (\log a_{i,i} - \log b_{i,i})^{2}}.
\]

Fix some nonprincipal ultrafilter $\mathcal{D}$ over $\omega$.
Let $X = P(m, \mathbb{R})$ and let $I \in X$ be the base point.
We next consider the question of which sequences 
of diagonal matrices lie in $X_{\infty}$. For each $n \geq 1$, let
\[
A_{n} = \Bigl( a^{(n)}_{i,j} \Bigr) \in X
\]
be a diagonal matrix. Then
\[
d_{n}(I,A_{n}) = \frac{ \sqrt{\sum (\log a^{(n)}_{i,i})^{2}} }{n}
\]
and so $(A_{n}) \in X_{\infty}$ if and only if there exists 
$k \geq 1$ such that
\begin{equation*} \label{E:exp}
e^{-kn} < a^{(n)}_{i,i} < e^{kn} \tag{Exp}
\end{equation*}
for each $1 \leq i \leq m$ and $n \geq 1$.

Before we can describe the structure of the asymptotic cone
$\Con_{\mathcal{D}}(X)$, we first need to recall the definition
of the corresponding field $\,^{*}\mathbb{R}_{\mathcal{D}}$ of 
nonstandard reals.

\begin{Def} \label{D:nonstandard}
Let $\mathbb{R}^{\omega}$ be the set of all sequences
$(x_{n})$ of real numbers. Define an equivalence relation
$\equiv$ on $\mathbb{R}^{\omega}$ by
\[
(x_{n}) \equiv (y_{n}) \quad \text{ if and only if } \quad
\{ n \in \omega \mid x_{n} = y_{n} \} \in \mathcal{D};
\]
and for each $(x_{n}) \in \mathbb{R}^{\omega}$,
let $[x_{n}]_{\mathcal{D}}$ be the corresponding equivalence
class. Then the field of {\em nonstandard reals\/} is defined
to be 
\[
\,^{*}\mathbb{R}_{\mathcal{D}} = 
\{ [x_{n}]_{\mathcal{D}} \mid (x_{n}) \in \mathbb{R}^{\omega} \}
\]
equipped with the operations
\begin{align*}
[x_{n}]_{\mathcal{D}} + [y_{n}]_{\mathcal{D}} &=
[x_{n} + y_{n}]_{\mathcal{D}} \\
[x_{n}]_{\mathcal{D}} \cdot [y_{n}]_{\mathcal{D}} &=
[x_{n} \cdot y_{n}]_{\mathcal{D}} 
\end{align*}
and the ordering
\[
[x_{n}]_{\mathcal{D}} < [y_{n}]_{\mathcal{D}} 
\quad \text{ if and only if } \quad
\{ n \in \omega \mid x_{n} < y_{n} \} \in \mathcal{D}.
\]
\end{Def}

In order to simplify notation, during the next few paragraphs,
we shall write $\,^{*}\mathbb{R}$ instead of the more precise
$\,^{*}\mathbb{R}_{\mathcal{D}}$. It is well-known that 
$\,^{*}\mathbb{R}$ is a nonarchimedean real closed field
and that $\mathbb{R}$ embeds into $\,^{*}\mathbb{R}$ via
the map $r \mapsto [r]_{\mathcal{D}}$, where 
$(r)$ denotes the sequence with constant value $r$.
(For the basic properties of $\,^{*}\mathbb{R}$, see
Lightstone-Robinson \cite{lr}.)

Now suppose that $(A_{n}) \in X_{\infty}$ is a sequence of diagonal 
matrices, where each $A_{n} = ( a^{(n)}_{i,j})$. Then we can define a
corresponding diagonal matrix
\[
(A_{n})_{\mathcal{D}} = \Bigl( \alpha_{i,j} \Bigr) \in 
P(m, \,^{*}\mathbb{R} )
\]
by setting $\alpha_{i,j} = [a^{(n)}_{i,j}]_{\mathcal{D}}$. Using
equation (\ref{E:exp}), we see that there exists $k \geq 1$ such that
each $\alpha_{i,i}$ satisfies the inequality
\[
\rho^{k} < \alpha_{i,i} < \rho^{-k},
\]
where $\rho \in \,^{*}\mathbb{R}$ is the positive
infinitesimal defined by
\[
\rho = ( e^{-n} )_{\mathcal{D}}.
\]
This suggests that $\Con_{\mathcal{D}}(X)$ should ``essentially''
be $P(m, K)$ for some field defined in terms of $\,^{*}\mathbb{R}$ 
and $\rho$.

\begin{Def} \label{D:rob}
Let $M_{0}$ be the subring of $\,^{*}\mathbb{R}$ defined by
\[
M_{0} = \{ t \in \,^{*}\mathbb{R} \mid
|t| < \rho^{-k} \text{ for some } k \geq 1 \}.
\]
Then $M_{0}$ has a unique maximal ideal
\[
M_{1} = \{ t \in \,^{*}\mathbb{R} \mid
|t| < \rho^{k} \text{ for all } k \geq 1 \}
\]
and the {\em Robinson field\/} $\,^{\rho}\mathbb{R}$
is defined to be the residue field $M_{0}/M_{1}$. 
By \cite{lr}, $\,^{\rho}\mathbb{R}$ is also a real closed field.
\end{Def}

If $(A_{n}) \in X_{\infty}$ is a sequence of diagonal matrices,
then
\[
(A_{n})_{\mathcal{D}} = \Bigl( \alpha_{i,j} \Bigr) \in 
P(m, M_{0}),
\]
and so we can define a corresponding matrix
\[
\overline{(A_{n})}_{\mathcal{D}} = 
\Bigl( \overline{\alpha}_{i,j} \Bigr) \in 
P(m, \,^{\rho}\mathbb{R} ),
\]
where each $\overline{\alpha}_{i,j} \in \,^{\rho}\mathbb{R}$ is
the element naturally associated with $\alpha_{i,j} \in \,^{*}\mathbb{R}$. 
Unfortunately, it is {\em not\/} always the case that if
$(A_{n})$, $(A^{\prime}_{n}) \in X_{\infty}$ correspond to the same
element of $\Con_{\mathcal{D}}(X)$, then 
$\overline{(A_{n})}_{\mathcal{D}} = \overline{(A^{\prime}_{n})}_{\mathcal{D}}$.
Thus, in order to obtain $\Con_{\mathcal{D}}(X)$, we must first factor
$P(m, \,^{\rho}\mathbb{R} )$ by a suitable equivalence relation.

\begin{Def} \label{D:val}
If $0 \neq \alpha \in M_{0}$, then $\log_{\rho} |\alpha|$
is a finite possibly nonstandard real and hence is infinitesimally close
to a unique standard real denoted by $\st ( \log_{\rho} |\alpha| )$.
By Lightstone-Robinson \cite[Section 3.3]{lr}, if 
$t \in M_{0} \smallsetminus M_{1}$ and $i \in M_{1}$,
then
\[
\st ( \log_{\rho} |t| ) = \st ( \log_{\rho} |t+i| ).
\]
Hence we can define a valuation
$\upsilon : \,^{\rho}\mathbb{R} \to \mathbb{R} \cup \{ \infty \}$ 
by
\[
\upsilon( \alpha ) = \st ( \log_{\rho} |\alpha| ).
\]
Let
\[
|\alpha|_{\upsilon} = e^{- \upsilon( \alpha )}
\]
be the associated absolute value.
\end{Def}

Following Leeb and Parreau \cite{p}, the asymptotic cone $\Con_{\mathcal{D}}(X)$
can now be described as follows. Let $V(m, \,^{\rho}\mathbb{R})$ be the
vector space of $m \times 1$ column vectors over the field 
$\,^{\rho}\mathbb{R}$. For each $A \in P(m,\,^{\rho}\mathbb{R})$, 
define the corresponding norm
\[
\varphi_{A} : V(m, \,^{\rho}\mathbb{R}) \to \mathbb{R}
\]
by
\[
\varphi_{A}(x) = 
\Big|\, \sqrt{ x^{t}A x}\;\Big|_{\upsilon}
\]
Then the map from $X_{\infty}$ defined by
$(A_{n}) \mapsto \varphi_{A}$, where 
$A = \overline{(A_{n})}_{\mathcal{D}}$, induces an isometry
\[
\Con_{\mathcal{D}}(X) \cong
\{ \varphi_{A} \mid A \in P(m,\,^{\rho}\mathbb{R}) \}.
\]
(When $A = ( \alpha_{i} ) \in P(m,\,^{\rho}\mathbb{R})$
and $B = ( \beta_{i} ) \in P(m,\,^{\rho}\mathbb{R})$ are
both diagonal matrices, then the corresponding distance in the
space of norms is given by
\[
d( \varphi_{A} , \varphi_{B} ) =
\sqrt{\sum ( \upsilon(\alpha_{i}) - \upsilon(\beta_{i} ))^{2}}.
\]
For more details, see Parreau \cite{p}.)
This space of norms is an instance of a classical construction of
Bruhat-Tits \cite{brt1}, \cite{brt2} and has a rich geometric structure. 
(The various
geometric notions discussed in the remainder of this paragraph will be 
defined and discussed in more detail in Section \ref{S:build}.) More 
precisely, $\Con_{\mathcal{D}}(X)$ is an affine $\mathbb{R}$-building; 
and the set $\mathcal{F}$ of apartments of $\Con_{\mathcal{D}}(X)$ is 
precisely the collection of subspaces of $\Con_{\mathcal{D}}(X)$ 
which are isometric to $\mathbb{R}^{m-1}$. Furthermore, the elements
of $\mathcal{F}$ correspond naturally to the unordered frames 
$\{ \,^{\rho}\mathbb{R}v_{1}, \cdots, \,^{\rho}\mathbb{R}v_{m} \}$
of the $m$-dimensional vector space $V(m, \,^{\rho}\mathbb{R})$
over $\,^{\rho}\mathbb{R}$. Let $\partial \Con_{\mathcal{D}}(X)$ be the 
associated spherical building at infinity of $\Con_{\mathcal{D}}(X)$.
Then the apartments of $\partial \Con_{\mathcal{D}}(X)$  are the boundaries at
infinity of the apartments of $\Con_{\mathcal{D}}(X)$; and so
the apartments of $\partial \Con_{\mathcal{D}}(X)$  also correspond 
naturally to the unordered frames 
$\{ \,^{\rho}\mathbb{R}v_{1}, \cdots, \,^{\rho}\mathbb{R}v_{m} \}$
of $V(m, \,^{\rho}\mathbb{R})$. As this observation suggests,
$\partial \Con_{\mathcal{D}}(X)$  is the usual spherical building 
associated with the natural $BN$-pair of $SL_{m}(\,^{\rho}\mathbb{R})$; 
i.e. $\partial \Con_{\mathcal{D}}(X)$  is the flag complex of 
the vector space $V(m, \,^{\rho}\mathbb{R})$. 
In particular, since $m \geq 3$, the isomorphism type of the Robinson
field $\,^{\rho}\mathbb{R}$ is determined by the isomorphism
type of the spherical building $\partial \Con_{\mathcal{D}}(X)$.
By the Kleiner-Leeb topological rigidity theorem for 
affine $\mathbb{R}$-buildings \cite{kl}, the isomorphism type of the 
spherical building $\partial \Con_{\mathcal{D}}(X)$  is a topological 
invariant of $\Con_{\mathcal{D}}(X)$. Consequently, the isomorphism type
of the Robinson field $\,^{\rho}\mathbb{R}$ is also a topological invariant
of $\Con_{\mathcal{D}}(X)$. 

From now on, we shall write  $\,^{\rho}\mathbb{R}_{\mathcal{D}}$ to
indicate the possible dependence on $\mathcal{D}$ of the Robinson field.
Since $\,^{\rho}\mathbb{R}_{\mathcal{D}}$ is a topological invariant
of $\Con_{\mathcal{D}}(X)$, the following result implies that 
if $CH$ fails, then $X = P(m, \mathbb{R})$ has $2^{2^{\omega}}$
asymptotic cones up to homeomorphism. Consequently, as the uniform lattice
$\Gamma$ is quasi-isometric to $P(m, \mathbb{R})$, it follows that
$\Gamma$ also has $2^{2^{\omega}}$ asymptotic cones up to homeomorphism. 

\begin{Thm} \label{T:rob}
If $CH$ fails, then there exists a set
$\{ \mathcal{D}_{\alpha} \mid \alpha < 2^{2^{\omega}} \}$ of
nonprincipal ultrafilters over $\omega$ such that
\[
\,^{\rho}\mathbb{R}_{\mathcal{D}_{\alpha}} \ncong
\,^{\rho}\mathbb{R}_{\mathcal{D}_{\beta}} 
\]
for all $\alpha < \beta < 2^{2^{\omega}}$. 
\end{Thm}

Now suppose that $CH$ holds. In this case, it is well-known that if
$\mathcal{A}$, $\mathcal{B}$ are nonprincipal ultrafilters, then
the corresponding fields $\,^{*}\mathbb{R}_{\mathcal{A}}$,
$\,^{*}\mathbb{R}_{\mathcal{B}}$ of nonstandard reals are
isomorphic. In \cite{t}, Thornton used the Diarra-Pestov \cite{d}, \cite{pe}
representation of $\,^{*}\mathbb{R}_{\mathcal{A}}$ as a Hahn
field to prove that if $CH$ holds, then the Robinson fields
$\,^{\rho}\mathbb{R}_{\mathcal{A}}$ and $\,^{\rho}\mathbb{R}_{\mathcal{B}}$
are isomorphic as valued fields. It follows that 
$X = P(m, \mathbb{R})$ has a unique asymptotic cone up to
isometry. In fact, Thornton proved that the analogous result holds for 
arbitrary symmetric spaces. Hence the following result holds for
uniform lattices in arbitrary connected semisimple Lie groups.
(Recall that the word metric $d$ on $\Gamma$ depends on the choice of
a finite generating set and so the metric space $( \Gamma, d)$
is determined only up to quasi-isometry. Consequently, the following
result is optimal.) 

\begin{Thm}[Thornton \cite{t}] 
Assume $CH$. If $\mathcal{A}$, $\mathcal{B}$
are nonprincipal ultrafilters, then
\[
\,^{\rho}\mathbb{R}_{\mathcal{A}} \cong \,^{\rho}\mathbb{R}_{\mathcal{B}}.
\]
Furthermore, if $G$ is a connected semisimple Lie group and  
$\Gamma$ is a uniform lattice in $G$, then $\Gamma$ has a unique
asymptotic cone up to bilipchitz homeomorphism.
\end{Thm}

Finally we should stress that it remains an open problem whether it
can be proved in $ZFC$ that there exists a finitely presented group 
with more than one asymptotic cone up to homeomorphism. However, 
when it comes to the question of whether there exists a finitely 
presented group with $2^{2^{\omega}}$ asymptotic cones up to 
homeomorphism, then the case is altered.

\begin{Thm} \label{T:ch}
If $CH$ holds, then every finitely generated group $\Gamma$ 
has at most $2^{\omega}$ asymptotic cones up to isometry.
\end{Thm}

\begin{Cor} \label{C:ch}
The following statements are equivalent.
\begin{enumerate}
\item[(a)] $CH$ fails.
\item[(b)] There exists a finitely presented group $\Gamma$
which has $2^{2^{\omega}}$ asymptotic cones up to
homeomorphism.
\item[(c)] There exists a finitely generated group $\Gamma$
which has $2^{2^{\omega}}$ asymptotic cones up to
homeomorphism.
\end{enumerate}
\end{Cor}

The rest of this paper is organised as follows. In Section
\ref{S:build}, we shall discuss the notions of an affine
$\mathbb{R}$-building and its spherical building at infinity;
and we shall explain why $\,^{\rho}\mathbb{R}_{\mathcal{D}}$ is a 
topological invariant of $\Con_{\mathcal{D}}(\Gamma)$, whenever
$\Gamma$ is a uniform lattice in a connected semisimple Lie group 
$G$ with at least one absolutely simple factor $S$ such that 
$\mathbb{R}\text{-rank}(S) \geq 2$. Sections \ref{S:code}
and \ref{S:ultra} will be devoted to the proof of
Theorem \ref{T:rob}. Finally we shall prove Theorem \ref{T:ch}
in Section \ref{S:ch}.

Throughout this paper, unless otherwise stated, the term ``Lie
group'' will always mean a real Lie group.

\section{Spherical and euclidean buildings} \label{S:build}

In this section, we shall discuss the notions of an affine
$\mathbb{R}$-building and its spherical building at infinity;
and we shall explain why $\,^{\rho}\mathbb{R}_{\mathcal{D}}$ is a 
topological invariant of $\Con_{\mathcal{D}}(\Gamma)$, whenever
$\Gamma$ is a uniform lattice in a connected semisimple Lie group 
$G$ with at least one absolutely simple factor $S$ such that 
$\mathbb{R}\text{-rank}(S) \geq 2$. 

Suppose that $V\cong\RR^n$ is a finite-dimensional Euclidean vector space,
with inner product $\langle-,-\rangle:V\times V\to\RR$.
For a nonzero vector $v\in V$, let $\sigma_v:V\to V$
denote the Euclidean reflection
$u\longmapsto u-2\frac{\langle u,v\rangle}{\langle v,v\rangle}v$.
A finite spanning set
$\Phi\subseteq V$ of (nonzero) vectors is called a root system
if $\sigma_uv\in\Phi$  and 
$2\frac{\langle x,v\rangle}{\langle v,v\rangle}\in\ZZ$ for
all $u,v\in\Phi$, see \cite[Ch.~VI]{bou}.
The group $W$ generated by the $\sigma_v$
is called the Weyl group of the root system.
Associated to such a finite reflection group $W$ is a certain
simplicial complex $\Sigma(W)$, its \emph{Coxeter complex}
\cite[Ch.~IV]{bou}:
roughly speaking, $\Sigma(W)$ is a $W$-invariant triangulation of the
unit sphere $\mathbb S^{n-1}\subseteq V$.
(The triangulation is obtained from the intersections
of $\mathbb S^{n-1}$ with the reflection hyperplanes
$v^\perp$ for $v\in\Phi$).

For example, the root system of type $A_n$ and its Coxeter complex
are defined as follows. Let
$\{e_1,\ldots,e_{n+1}\}$ denote the standard orthonormal basis for 
$\RR^{n+1}$ and let $V$ be the subspace of $\RR^{n+1}$ defined by
\[
V=\{x\in\RR^{n+1}\mid x_1+\cdots +x_{n+1}=0\}.
\]
Then $\Phi=\{e_i-e_j\mid i\neq j\} \subseteq V$ is a root system
of type $A_{n}$. As an abstract group, the Weyl group
$W$ is the symmetric group on $n+1$ letters, acting by coordinate
permutations; and as a poset,
$\Sigma(W)$ is the set of all
$\subseteq$-ordered chains consisting
of nontrivial subsets of $\{1,\ldots,n+1\}$
(i.e. $\Sigma(W)$ is the first barycentric subdivision of
the boundary of an $n+1$-simplex).

A \emph{spherical building} is an abstract simplicial complex 
(a poset) $(\Delta,\leq)$
with a distinguished collection of subcomplexes $\Sigma$, called
\emph{apartments}, satisfying the following axioms.\\
\textbf{(B$_1$)} Each apartment $\Sigma\subseteq\Delta$
is isomorphic to a (fixed) Coxeter complex $\Sigma(W)$.\\
\textbf{(B$_2$)} Any two simplices in $\Delta$
are contained in some apartment.\\
\textbf{(B$_3$)} Given two apartments $\Sigma_1,\Sigma_2\subseteq\Delta$,
there exists an isomorphism $\Sigma_1\to\Sigma_2$ fixing
$\Sigma_1\cap\Sigma_2$ simplex-wise.\\
For more details, we refer to \cite{br}, \cite{ro}, \cite{ti}.
The standard examples of buildings are obtained from algebraic
groups as follows.

Let $\uG$ be a reductive algebraic group defined over a field $k$
\cite{b} \cite{bt} \cite{sp}
and let $\Delta$ denote the poset of all $k$-parabolic subgroups
of $\uG$, ordered by reversed inclusion. Then
$\Delta=\Delta(\uG,k)$ is a spherical building, the canonical
building associated to $\uG$ over $k$, and the group $\uG(k)$ 
of $k$-points of $\uG$ acts strongly transitively on $\Delta(\uG,k)$,
see Tits \cite[Ch.~5]{ti}.
If $\uG$ is absolutely simple (i.e. if $\uG$ is simple over the
algebraic closure $\bar k$) and adjoint, then
the building uniquely determines the field of definition $k$.
\begin{Thm}
\cite[5.8]{ti}
\label{FieldsAreEqual}
Let $\uG$ and $\uG'$ be adjoint absolutely simple algebraic groups 
of rank at least $2$ defined over the fields $k$ and $k'$. 
If there is a building isomorphism $\Delta(\uG,k)\cong\Delta(\uG',k')$,
then the fields $k$ and $k'$ are isomorphic.
\qed
\end{Thm}
We next introduce the notion of an affine $\RR$-building (also called a
\emph{Euclidean building}). Below we shall give Tits' definition \cite{ti2},
as corrected in \cite[App.~3]{ro}.
(A different set of axioms was proposed in Kleiner-Leeb \cite{kl},
based on nonpositive curvature and geodesics; in \cite{p},
Parreau showed that these two approaches are equivalent).
Let $W$ be the Weyl group of a root system $\Phi\subseteq V$,
and let $W_\aff$ denote the
semidirect product of $W$ and the vector group $(V,+)$.
Then in $V$, we obtain the corresponding \emph{reflection hyperplanes}
(the fixed point sets of reflections),
\emph{half-spaces} (determined by reflection hyperplanes), and 
\emph{Weyl chambers} (the fundamental domains for $W$), see \cite{bou}.
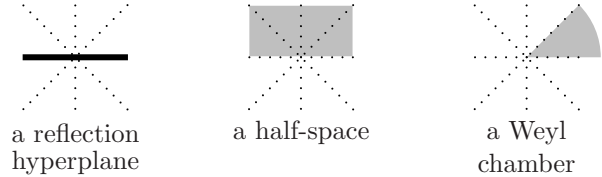
\begin{center}
\begin{pspicture}(8,2)
\psframe[linestyle=none,fillstyle=solid,fillcolor=lightgray](3.3,1)(4.7,1.7)
\pswedge[linestyle=none,fillstyle=solid,fillcolor=lightgray](7,1){1}{0}{45}
\psline[linestyle=dotted](.3,.3)(1.7,1.7)
\psline[linewidth=0.8mm](.3,1)(1.7,1)
\psline[linestyle=dotted](1,.3)(1,1.7)
\psline[linestyle=dotted](.3,1.7)(1.7,.3)
\psline[linestyle=dotted](3.3,.3)(4.7,1.7)
\psline[linestyle=dotted](3.3,1)(4.7,1)
\psline[linestyle=dotted](4,.3)(4,1.7)
\psline[linestyle=dotted](3.3,1.7)(4.7,.3)
\psline[linestyle=dotted](6.3,.3)(7.7,1.7)
\psline[linestyle=dotted](6.3,1)(7.7,1)
\psline[linestyle=dotted](7,.3)(7,1.7)
\psline[linestyle=dotted](6.3,1.7)(7.7,.3)
\rput(1,0){a reflection}
\rput(1,-.4){hyperplane}
\rput(4,0){a half-space}
\rput(7,0){a Weyl}
\rput(7,-.4){chamber}
\end{pspicture}\\\
\end{center}
\begin{Def}
Fix $W$ and $V$ as above.
A pair $(\cI,\cF)$ consisting of a nonempty set $\cI$ and a family $\cF$ of
injections $\phi:V\to\cI$ is called an
\emph{affine $\RR$-building} if it has the following properties.\\
\textbf{(A$\RR$B$_1$)}
If $w\in W_\aff$ and $\phi\in\cF$, then $\phi\circ w\in\cF$.\\
\textbf{(A$\RR$B$_2$)}
For $\phi,\psi\in\cF$, the preimage $X=\phi^\inv\psi(V)$ is closed
and convex (possibly empty), and there exists $w\in W_\aff$ 
such that $\phi$ and $\psi\circ w$ agree on $X$.\\
\textbf{(A$\RR$B$_3$)}
Given $x,y\in\cI$, there exists $\phi\in\cF$ with $\{x,y\}\subseteq\phi(V)$.\\
Any $\phi$-image of $V$ is called an \emph{apartment}; the $\phi$-image of
a reflection hyperplane, a half-space, and a Weyl chamber is called a
\emph{wall}, a \emph{half-apartment}, and a \emph{sector}, respectively.
A wall is \emph{thick} if it bounds three distinct half-apartments, and
a point is \emph{thick} if every wall passing through it is thick.\\
\textbf{(A$\RR$B$_4$)}
Given two sectors $S_1,S_2\subseteq\cI$, there exist subsectors
$S_i'\subseteq S_i$ and an apartment $F\subseteq\cI$ with
$S_1'\cup S_2'\subseteq F$.\\
\textbf{(A$\RR$B$_5$)}
If $F_1,F_2,F_3$ are apartments having pairwise a half-apartment in
common, then $F_1\cap F_2\cap F_3\neq\emptyset$.\\
The \emph{dimension} of an affine $\RR$-building is the vector space
dimension of $V$.
\end{Def}
It follows that $\cI$ admits a unique metric $d$ which pulls back to the
Euclidean metric on $V$ for every $\phi:V\to\cI$, and
that $(\cI,d)$ is a CAT$(0)$-space. (See \cite{bh} for the notion
of a CAT$(0)$-space.)
We say that $\cI$ is \emph{complete} if
$(\cI,d)$ is complete as a metric space (every Cauchy sequence converges)
and that $\cF$ is \emph{complete} if every injection $\phi:V\to\cI$
which is compatible with the axioms (A$\RR$B$_1$)--(A$\RR$B$_5$) is
already in $\cF$. It can be shown that every $\cF$ admits a unique
completion $\widehat\cF$ such that $(\cI,\widehat\cF)$ is
an affine $\RR$-building;
the metric completion of $\cI$, however, is in general not
an affine $\RR$-building.
\begin{Prop}
\label{HomeoPresFlats}
Suppose that $(\cI,\cF)$ and $(\cI',\cF')$ are affine $\RR$-buildings, and that
$f:\mathcal I\to\mathcal I'$ is a homeomorphism. Let $F\subseteq\cI$
be an apartment. If $\cI'$ and $\cF'$ are complete, then $f(F)$ is an
apartment in $\mathcal I'$.
\end{Prop}
\begin{proof}
This was first proved by Kleiner-Leeb \cite[Prop. 6.4.1]{kl}.  
A sheaf-theoretic proof was given in Kramer-Tent \cite{kt}, 
using the fact that apartments can be viewed as certain
global sections in the orientation sheaf of
the topological space $\mathcal I$, which can be characterised topologically.
From this, one deduces that both buildings have the same dimension,
and that
$f(F)$ is locally isometric to Euclidean space.
Since $\cI'$ is complete,  $f(F)$ is also complete
and thus -- being contractible -- globally isometric to Euclidean space.
Such a subspace is always an apartment in the completion of $\cF'$.
\end{proof}
Proposition \ref{HomeoPresFlats} and Proposition \ref{P:drafts}
below are false if $\cI'$ is not assumed to be complete.

An affine $\RR$-building $(\cI,\cF)$ has a
\emph{spherical building at infinity},
denoted by $\partial(\cI,\cF)$. The chambers of this building
are equivalence classes of sectors, where two sectors are equivalent if
the Hausdorff distance between them is finite. If a basepoint $x_0\in\cI$ is
fixed, then every chamber of $\partial(\cI,\cF)$ has a unique
$x_0$-based sector as its representative.

The \emph{draft} $D(\cI,\cF)$ of an affine $\RR$-building $(\cI,\cF)$ is
the pair $(\cI,\cA)$ consisting of all points and all apartments
of $(\cI,\cF)$.
As a direct consequence of \ref{HomeoPresFlats}, we have the following
result.
\begin{Prop} \label{P:drafts}
Let $(\cI,\cF)$ and $(\cI',\cF')$ be affine $\RR$-buildings such that
$\cI$, $\cI'$, $\cF$ and $\cF'$ are complete. Then
each homeomorphism $f:\cI\to\cI'$
induces an isomorphism of drafts $D(\cI,\cF)\rTo^\cong D(\cI',\cF')$.
\qed
\end{Prop}
The draft disregards the metric structure and
the Weyl group of the affine $\RR$-building. In general, nonisomorphic
affine $\RR$-buildings can have isomorphic drafts. However, we have the
following result.
\begin{Prop}
\label{DraftsAtInfinity}
Suppose that $\cI$ contains a thick point $x$.
If 
\[f:D(\cI,\cF)\rTo^\cong D(\cI',\cF') 
\]
is an isomorphism
of drafts of affine $\RR$-buildings,
then $f$ induces an isomorphism
\[
\partial f:\partial(\cI,\cF)\to\partial(\cI',\cF') 
\]
between the
respective spherical buildings at infinity.
\end{Prop}
\begin{proof}
Fix some thick point $x$ of $\cI$.
For each $y\in\mathcal I$, let $\cvx\{x,y\}$ denote the intersection of all
apartments containing $x$ and $y$. This set should be
pictured as a diamond-shaped set with $x$ as one tip.
Since $x$ has thick walls, $\cvx\{x,y\}$ is always contained
in one of the $x$-based sectors of $\cI$.
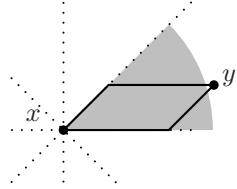
\begin{center}
\begin{pspicture}(3,3)
\pswedge[linestyle=none,fillstyle=solid,fillcolor=lightgray](1,1){2}{0}{45}
\pspolygon(1,1)(2.4,1)(3,1.6)(1.6,1.6)
\psline[linestyle=dotted](.3,.3)(2.7,2.7)
\psline[linestyle=dotted](.3,1)(2.7,1)
\psline[linestyle=dotted](1,.3)(1,2.7)
\psline[linestyle=dotted](.3,1.7)(1.7,.3)
\psdots(1,1)(3,1.6)
\rput(.6,1.2){$x$}
\rput(3.2,1.7){$y$}
\end{pspicture}
\end{center}
Let $F\subseteq\cI$ be an apartment containing $x$. In the poset
$\left(\left\{\cvx\{x,y\}\mid y\in F\right\},\subseteq\right)$,
the unions over the maximal chains are precisely the $x$-based sectors
in $F$. Thus any isomorphism of drafts preserves $x$-based sectors.
The intersections of the $x$-based sectors in $F$ naturally form 
a poset isomorphic to $\Sigma(W)$; and the union of these posets,
where $F$ runs through all apartments containing $x$, is canonically
isomorphic to the spherical building $\partial(\cI,\cF)$.
The result now follows from the fact that the underlying poset of a spherical 
building completely determines the building itself.
\end{proof}

Let $\uG$ be a semisimple algebraic group defined over the real closure
$\QR$ of $\QQ$. Then the \emph{$\mathbb{R}$-rank} of $\uG$ is the maximal
dimension of a $\QR$-split torus $\uS$ of $\uG$.
If $\cR$ is any real closed field (such as $\cR=\RR$ or
$\cR={}^\rho\RR$),
then $\QR\subseteq\cR$ and $\uS$ is also a maximal $\cR$-split torus over
$\cR$. In fact, $\uG$ has the same structure,
rank, Tits diagram and building type over $\cR$ as over its field
of definition $\QR$, see \cite[Sec.~4]{gr}. Note that the group
$G=\uG(\RR)$ of $\RR$-points of $\uG$, endowed with the Hausdorff topology,
is a real Lie group; and its (Hausdorff) connected component $G^\circ$ is
a semisimple Lie group $G$ such that $[G : G^\circ] < \infty$. 
Furthermore, if $G$ is any connected semisimple Lie group, then
there exists a semisimple algebraic group defined over $\QR$ such
that $G/Z(G)$ and $\uG(\RR)^\circ$ are isomorphic as Lie groups; for
example, see ~\cite[1.14.6]{eb}. (In fact, $\uG$ can even be taken to be
defined over $\QQ$; but in our setting, it is more convenient to work with
real closed fields.) We then define $\mathbb{R} \text{-rank}(G)$ to be
the $\mathbb{R}$-rank of the algebraic group $\uG$; and we define $G$
to be absolutely simple if and only if $\uG$ is absolutely simple.
(Equivalently, $G$ is absolutely simple if and only if the
complexification $\mathfrak{g} \otimes_{\mathbb{R}} \mathbb{C}$
of the Lie algebra $\mathfrak{g}$ of $G$ remains simple.)
Since $Z(G)$ is finite, the Riemannian manifolds $G$ and $G/Z(G)$ are
quasi-isometric. Hence, in the remainder of this section, we can
restrict our attention to connected semisimple Lie groups of the
form $\uG(\RR)^\circ$.

Fix a semisimple algebraic group $\uG$ defined over $\QR$
of $\mathbb{R}$-rank $m\geq 1$.
Let $\uS\subseteq\uG$ be a maximal $\QR$-split torus
and $\uN=\mathrm{Nor}_\uG(\uS)$ its normaliser. The quotient
$W=\uN/\uS$ is the relative Weyl group for $\uG$ 
(once again, over any real closed field $\cR$). Let
$\cD$ be a non-principal ultrafilter, let $\eps=(1/n)_\cD$, and let
$\rho=e^\eps=(e^{-n})_\cD$.
The Robinson field $^\rho\RR$ has a unique maximal
o-convex subring $^\rho\cO\subseteq{}^\rho\RR$,
corresponding to the canonical o-valuation
$\nu:{}^\rho\RR\to\RR\cup\{\infty\}$, which was defined in 
Definition \ref{D:val}. (For general results about
real closed fields and o-valuations, see \cite{pr}.)
Since $^\rho\cO$ is o-convex, we have that
$\QR\subseteq{}^\rho\cO\subseteq{}^\rho\RR$.
Therefore the groups $\uG({}^\rho\cO)\subseteq\uG({}^\rho\RR)$ and
$\uS({}^\rho\cO)\subseteq\uS({}^\rho\RR)$ are defined, as well as
the coset spaces
\[
\cI=\uG({}^\rho\mathbb R)/\uG({}^\rho\cO)\quad\text{ and }\quad
V=\uS({}^\rho\mathbb R)/\uS({}^\rho\cO).
\]
Moreover, $V$ can be regarded in a natural way as an $m$-dimensional
real vector space, equipped with a natural action of $W$ as a finite 
reflection group. 
For example, if $\uG=SL_{m+1}$ and $\uS$ is the group of diagonal
matrices in $SL_{m+1}$, then $\uN$ consists of permutation
matrices acting by coordinate permutations; and the resulting root system
is the one of type $A_m$ described at the beginning of this section.

If we extend the $W$-action by the translations, then we obtain an
affine Weyl group $W_\aff$. For $g\in\uG({}^\rho\RR)$
and $v=n\uS({}^\rho\cO)\in\uS({}^\rho\RR)/\uS({}^\rho\cO)=V$, let 
$\phi_g(v)=gn\uG({}^\rho\cO)\in\cI$ and put
\[
\cF=\{\phi_g\mid g\in\uG({}^\rho\mathbb R_\cD)\}.
\]
\begin{Thm}
The pair $(\cI,\cF)=\Delta_\aff(\uG,{}^\rho\mathbb R,{}^\rho\cO)$ is
an affine $\RR$-building such that $\cI$ and $\cF$ are complete.
The building at infinity is the spherical building
$\partial\Delta_\aff(\uG,{}^\rho\mathbb R,{}^\rho\cO)=
\Delta(\uG,{}^\rho\mathbb R)$.
\end{Thm}
\begin{proof}
This is a special case of a much
more general result on the affine $\Lambda$-buildings
associated to arbitrary real closed valued fields, which was 
proved in Kramer-Tent \cite{kt}.
The fact that $\cI$ and $\cF$ are complete follows from the
$\omega_1$-saturatedness of countable ultrapowers.
\end{proof}
Now we shall explain how asymptotic cones fit into the picture.
Let $G$ be a semisimple Lie group of rank $m\geq 1$,
endowed with a left-invariant Riemannian metric $d$.
Let $K \leqslant G$ be a maximal compact subgroup. The Riemannian
symmetric space $X=G/K$ carries a natural metric (unique up to
homothety), and it is not difficult to show that the natural map
$G\to G/K$ is a quasi-isometry. Thus $X$ and $G$ have
bilipschitz homeomorphic asymptotic cones.
As above, let $\uG$ be an algebraic group over
$\QR$ with $\uG(\RR)^\circ=G$.

\begin{Prop} \label{P:affine}
The asymptotic cone $\Con_\cD(X)$ is isometric to the point space
$\cI$ of the building $\Delta_\aff(\uG,{}^\rho\RR_\cD,{}^\rho\cO_\cD)$.
\end{Prop}
\begin{proof}
This was proved in Kramer-Tent \cite{kt} and also independently
in Thornton \cite{t}.
\end{proof}

The fact that $\Con_\cD(X)$ is an affine
$\RR$-building was proved first by Kleiner-Leeb \cite{kl}.
(A vague conjecture pointing in this direction was made by Gromov in
\cite[p.~54]{g2}). However, Kleiner and Leeb
did not determine the building which one obtains. 
As we mentioned above, the fact that
$\Con_\cD(X)$ can be identified as a metric space with the quotient
$\uG({}^\rho\mathbb R)/\uG({}^\rho\cO)$
was proved independently by Thornton \cite{t}, but he did not
identify the affine building or the spherical building at infinity.
This was done by Leeb and Parreau \cite{p} for the special case of
$\uG=SL_{m+1}$; and
one should also mention Bennett's general result \cite{be} on the
affine $\Lambda$-buildings related to $SL_{m+1}(F)$ for arbitrary
valued fields $F$. 

\begin{Thm} \label{ConeAndField}
Suppose that $G$ is a connected absolutely simple Lie group such that
$\mathbb{R}\text{-rank}(G) \geq 2$ and let $\Gamma$ be a
uniform lattice in $G$. Let $\cD$ and $\cD'$ be nonprincipal ultrafilters.
If $\Con_\cD(\Gamma)$ and $\Con_{\cD'}(\Gamma)$ are homeomorphic, then
the Robinson fields $^\rho\RR_\cD$ and $^\rho\RR_{\cD'}$ are isomorphic.
\end{Thm}
\begin{proof}
Let $\uG$ be an absolutely simple algebraic group defined over
$\QR$ such that $\mathbb{R}\text{-rank}(\uG) \geq 2$ and
$G = \uG(\mathbb{R})^{\circ}$. Let $K \leqslant G$ be a maximal 
compact subgroup. Since $\Gamma$ is quasi-isometric to 
the symmetric space $X = G/K$, the asymptotic cones
$\Con_\cD(\Gamma)$ and $\Con_{\cD'}(\Gamma)$ are bilipschitz
homeomorphic to $\Con_\cD(X)$, $\Con_{\cD'}(X)$ respectively.
Consequently, if $\Con_\cD(\Gamma)$ and $\Con_{\cD'}(\Gamma)$ are
homeomorphic, then $\Con_\cD(X)$ and $\Con_{\cD'}(X)$ are also
homeomorphic. Hence, by Propositions \ref{P:affine} and \ref{P:drafts}, 
the affine $\RR$-buildings
$\Delta_\aff(\uG,{}^\rho\RR_\cD,{}^\rho\cO_\cD)$ and
$\Delta_\aff(\uG,{}^\rho\RR_{\cD'},{}^\rho\cO_{\cD'})$ have isomorphic drafts.
Applying Proposition \ref{DraftsAtInfinity}, it follows that
the corresponding buildings at infinity 
$\Delta(\uG,{}^\rho\RR_\cD)$ and $\Delta(\uG,{}^\rho\RR_{\cD'})$ are
also isomorphic.
Finally, by Theorem \ref{FieldsAreEqual}, the building
$\Delta(\uG,{}^\rho\RR_\cD)$
determines the field $^\rho\RR_\cD$ up to isomorphism, and so
${}^\rho\RR_\cD\cong{}^\rho\RR_{\cD'}$ .
\end{proof}
We should make a few comments concerning the hypotheses
on the Lie group $G$ in the statement of Theorem \ref{ConeAndField}.
If $G$ has $\mathbb{R}$-rank $1$, then $\Con_\cD\Gamma$
is a homogeneous $\RR$-tree with uncountable branching at every
point. Furthermore, the isometry type of this $\mathbb{R}$-tree
is independent of the choice of the ultrafilter $\cD$ (and even
of the Lie type of $G$.) For example, see Dyubina-Polterovich \cite{dp}. 
Thus the hypothesis on the $\mathbb{R}$-rank of $G$ is
certainly necessary. To understand the reason for the hypothesis
that $G$ is absolutely simple, consider the complex Lie group
$G = SL_{n}(\mathbb{C})$ for some $n \geq 3$. By embedding
$SL_{n}(\mathbb{C})$ as an algebraic subgroup 
$\uG(\mathbb{R})$ of $GL_{2n}(\mathbb{R})$, we can regard $G$
as a real simple Lie group of $\mathbb{R}$-rank $n-1$.
However, $G$ is not absolutely simple, since 
$\uG(\mathbb{C}) \cong SL_{n}(\mathbb{C}) \times SL_{n}(\mathbb{C})$.
(More generally, it turns out that a real simple Lie group $G$
is absolutely simple if and only if $G/Z(G)$ is not isomorphic to a 
complex Lie group.) When we consider the spherical building at infinity
of the corresponding affine $\mathbb{R}$-building, then we
are only able to recover the algebraic
closure $^\rho\mathbb{C}_{\cD}={}^\rho\RR_{\cD}(\sqrt{-1})$ 
rather than the Robinson field $^\rho\RR_\cD$ itself.
Since $^\rho\mathbb{C}_{\cD}$ is an algebraically closed
field of cardinality $2^{\omega}$, it follows that
$^\rho\mathbb{C}_{\cD} \cong \mathbb{C}$ for every
nonprincipal ultrafilter $\cD$ over $\omega$. Thus the
fields $^\rho\mathbb{C}_{\cD}$ cannot be used to distinguish
between the (possibly different) asymptotic cones
of $SL_{n}(\mathbb{C})$; and  it remains an open question whether 
the asymptotic cones of complex simple 
Lie groups depend on the chosen ultrafilter. However,
there is an obvious generalisation of Theorem \ref{ConeAndField}
to semisimple Lie groups.
\begin{Cor}
Let $G$ be a connected semisimple Lie group with at least
one absolutely simple factor $S$ such that 
$\mathbb{R}\text{-rank}(S) \geq 2$ and let $\Gamma$ be a
uniform lattice in $G$. Let $\cD$ and $\cD'$ be nonprincipal 
ultrafilters. If
$\Con_\cD(\Gamma)$ and $\Con_{\cD'}(\Gamma)$ are homeomorphic, then
the corresponding Robinson fields $^\rho\RR_\cD$ and $^\rho\RR_{\cD'}$ 
are isomorphic.
\end{Cor}
This follows from Theorem \ref{ConeAndField}, together with
the fact that the buildings at infinity of $\Con_\cD(\Gamma)$ 
and $\Con_{\cD'}(\Gamma)$ decompose into products of the buildings 
corresponding to the simple factors of $G$.

\section{Invariants of linear orders} \label{S:code}

The next two sections will be devoted to the proof of 
Theorem \ref{T:rob}. We will begin by reducing Theorem \ref{T:rob}
to an analogous statement concerning the linearly ordered sets
$\omega^{\omega}/\mathcal{D}$, where $\mathcal{D}$ is a
nonprincipal ultrafilter over $\omega$. 

\begin{Def} \label{D:filter}
A {\em filter\/} over $\omega$ is a collection
$D$ of subsets of $\omega$ satisfying the
following conditions:
\begin{enumerate}
\item[(i)] $\omega \in D$.
\item[(ii)] If $A$, $B \in D$, then 
$A \cap B \in D$.
\item[(iii)] If $A \in D$ and 
$A \subseteq B \subseteq \omega$, then
$B \in D$.
\end{enumerate}
The filter $D$ is {\em nontrivial\/} if and only if
\begin{enumerate}
\item[(iv)] $\emptyset \notin D$.
\end{enumerate}
\end{Def}

Let $D$ be a nontrivial filter over $\omega$. Then $\equiv_{D}$
is the equivalence relation defined on 
$\omega^{\omega} = \{ f \mid f: \omega \to \omega \}$ by
\[
f \equiv_{D} g \quad \text{ if and only if } \quad 
\{ n \in \omega \mid f(n) = g(n) \} \in D.
\]
For each $f \in \omega^{\omega}$, we denote the corresponding
$\equiv_{D}$-equivalence class by $f/D$; and we let
\[
\omega^{\omega}/D = \{ f/D \mid f \in \omega^{\omega} \}
\]
equipped with the partial order defined by
\[
f/D < g/D \quad \text{ if and only if } \quad 
\{ n \in \omega \mid f(n) < g(n) \} \in D.
\]
As usual, we identify each natural number $\ell \in \omega$ with the
corresponding element $c_{\ell}/D \in \omega^{\omega}$, defined by
$c_{\ell}(n) = \ell$ for all $n \in \omega$. If $\mathcal{D}$ is a
nonprincipal ultrafilter over $\omega$, then 
$\omega^{\omega}/\mathcal{D}$ is a linear order. In this case,
we define
\[
( \omega^{\omega}/ \mathcal{D} )^{*} = \{ g/\mathcal{D} \in 
\omega^{\omega}/ \mathcal{D} \mid \ell < g/\mathcal{D} 
\text{ for all } \ell \in \omega \}.
\]

As we shall now explain, Theorem \ref{T:rob} is an easy consequence of
Theorem \ref{T:point}, which we shall prove in Section \ref{S:ultra}.

\begin{Def} \label{D:iso}
Let $L_{1}$, $L_{2}$ be linear orders.
\begin{enumerate}
\item[(a)] $L_{1} \approx_{f} L_{2}$ if and only if $L_{1}$ and $L_{2}$ have nonempty
isomorphic final segments. 
\item[(b)] $L_{1} \approx_{i} L_{2}$ if and only if $L_{1}$ and $L_{2}$ have nonempty
isomorphic initial segments. 
\end{enumerate}
\end{Def}

\begin{Thm} \label{T:point}
If $CH$ fails, then there exists a set
$\{ \mathcal{D}_{\alpha} \mid \alpha < 2^{2^{\omega}} \}$ of
nonprincipal ultrafilters over $\omega$ such that
\[
( \omega^{\omega}/ \mathcal{D}_{\alpha} )^{*} \not \approx_{i}
( \omega^{\omega}/ \mathcal{D}_{\beta} )^{*} 
\]
for all $\alpha < \beta < 2^{2^{\omega}}$.
\end{Thm}

\begin{proof}[Proof of Theorem \ref{T:rob}]
For each nonprincipal ultrafilter $\mathcal{D}$ over $\omega$, 
let
\[
\,^{\rho}\mathbb{R}^{\infty}_{\mathcal{D}} =
\{ a \in \,^{\rho}\mathbb{R}_{\mathcal{D}} \mid a > \ell 
\text{ for every } \ell \in \omega \};
\]
and let $E_{\mathcal{D}}$ be the convex equivalence relation
defined on $\,^{\rho}\mathbb{R}^{\infty}_{\mathcal{D}}$ by
\[
a \mathrel{E_{\mathcal{D}}} b \quad \text{ if and only if } \quad
|a - b| < \ell \text{ for some } \ell \in \omega.
\]
Let $L_{\mathcal{D}} = 
\,^{\rho}\mathbb{R}^{\infty}_{\mathcal{D}}/E_{\mathcal{D}}$,
equipped with the quotient linear ordering; and regard 
$L_{\mathcal{D}} \times \mathbb{Z}$ as a linear ordering with
respect to the usual lexicographical ordering, defined by
$(a_{1},z_{1}) < (a_{2},z_{2})$
if and only if either:
\begin{itemize}
\item $a_{1} < a_{2}$; or
\item $a_{1} = a_{2}$ and $z_{1} < z_{2}$.
\end{itemize}
Then it is easily checked that
\[
L_{\mathcal{D}} \times \mathbb{Z} \cong 
\{ g/\mathcal{D} \in ( \omega^{\omega}/ \mathcal{D} )^{*} \mid
g/\mathcal{D} < \rho^{-n} \text{ for some } n \geq 1 \}.
\]
Now suppose that $\mathcal{A}$, $\mathcal{B}$ are nonprincipal 
ultrafilters over $\omega$ and that 
$f: \,^{\rho}\mathbb{R}_{\mathcal{A}} \to \,^{\rho}\mathbb{R}_{\mathcal{B}}$
is a field isomorphism. Since $\,^{\rho}\mathbb{R}_{\mathcal{A}}$,
$\,^{\rho}\mathbb{R}_{\mathcal{B}}$ are real closed, it follows that
$f$ is also order-preserving. It is also clear that
$f [ \,^{\rho}\mathbb{R}^{\infty}_{\mathcal{A}} ] =
\,^{\rho}\mathbb{R}^{\infty}_{\mathcal{B}}$ and that $f$ maps the
equivalence relation $E_{\mathcal{A}}$ to the equivalence relation 
$E_{\mathcal{B}}$. It
follows that the ordered sets $L_{\mathcal{A}}$ and $L_{\mathcal{B}}$
are isomorphic; and hence that
\[
( \omega^{\omega}/ \mathcal{A} )^{*} \approx_{i}
( \omega^{\omega}/ \mathcal{B} )^{*}. 
\]
Consequently, if $\{ \mathcal{D}_{\alpha} \mid \alpha < 2^{\kappa} \}$ 
is the set of nonprincipal ultrafilters over $\omega$ given by
Theorem \ref{T:point}, then 
\[
\,^{\rho}\mathbb{R}_{\mathcal{D}_{\alpha}} \ncong
\,^{\rho}\mathbb{R}_{\mathcal{D}_{\beta}} 
\]
for all $\alpha < \beta < 2^{\kappa}$. 
\end{proof}

A similar argument shows that the corresponding fields
$\,^{*}\mathbb{R}_{\mathcal{D}_{\alpha}}$, 
$\alpha < 2^{\kappa}$, of nonstandard reals
are also pairwise nonisomorphic.
This improves a result of Roitman \cite{r}, who proved that
it is consistent that there exist $2^{\omega}$ pairwise
nonisomorphic fields of nonstandard reals, each of the
form $\,^{*}\mathbb{R}_{\mathcal{D}}$ for some
nonprincipal ultrafilter $\mathcal{D}$ over $\omega$. 

Most of the remainder of this section will be devoted to
the construction of a collection of extremely nonisomorphic
linear orders, which will later be used as suitable ``invariants''
in the proof of Theorem \ref{T:point}. This construction is
a special case of the more general techniques which are
developed in Chapter III of Shelah \cite{s2}. In order to make
this paper relatively self-contained, we have provided 
proofs of the relevant results. We assume that the reader is familiar
with the basic properties of regular cardinals, singular cardinals, and
stationary subsets of regular cardinals. (For example, see Sections 6 
and 7 of Jech \cite{j}.)

\begin{Def} \label{D:cf}
Suppose that $I$ is a linear order and that 
$\emptyset \neq A \subseteq I$.
\begin{enumerate}
\item[(a)] A subset $B \subseteq A$ is said to be {\em cofinal\/}
in $A$ if for all $a \in A$, there exists an element $b \in B$
such that $a \leq b$. The {\em cofinality\/} of $A$ is defined to
be 
\[
\cf(A) = \min \{ |B| \mid B \text{ is a cofinal subset of } A \}.
\]
\item[(b)] A subset $B \subseteq A$ is said to be {\em coinitial\/}
in $A$ if for all $a \in A$, there exists an element $b \in B$
such that $b \leq a$. The {\em coinitiality\/} of $A$ is defined to
be 
\[
\ci(A) = \min \{ |B| \mid B \text{ is a coinitial subset of } A \}.
\]
\end{enumerate}
\end{Def}

\begin{Def} \label{D:cut}
Suppose that $I$ is a linear order and that $\lambda$, $\theta \geq \omega$
are regular cardinals. Then $(I_{1}, I_{2})$ is a $(\lambda, \theta)$-cut 
of $I$ if the following conditions hold:
\begin{enumerate}
\item[(a)] $I = I_{1} \cup I_{2}$ and $s < t$ for all
$s \in I_{1}$, $t \in I_{2}$.
\item[(b)] $\cf(I_{1}) = \lambda$.
\item[(c)] $\ci(I_{2}) = \theta$. 
\end{enumerate}
\end{Def}

\begin{Def} \label{D:inv}
Suppose that $J$, $L$ are linear orders. Then the order-preserving map
$\varphi: J \to L$ is an {\em invariant embedding\/} if whenever
$(J_{1}, J_{2})$ is a $(\lambda, \theta)$-cut of $J$ for some
$\lambda$, $\theta > \omega$, then there does {\em not\/}
exist an element $x \in L$ such that $\varphi(s) < x < \varphi(t)$ for all
$s \in J_{1}$, $t \in J_{2}$.

In this case, $\varphi$ is said to be an {\em invariant cofinal embedding\/}
if $\varphi[J]$ is cofinal in $L$ and $\varphi$ is said to be an 
{\em invariant coinitial embedding\/} if $\varphi[J]$ is coinitial in $L$ 
\end{Def}

\begin{Lem} \label{L:code}
If $\lambda > \omega_{1}$ is a regular cardinal, then there exists a set 
$\{ I_{\alpha} \mid \alpha < 2^{\lambda} \}$ of linear orders satisfying the 
following conditions:
\begin{enumerate}
\item[(a)] $\cf(I_{\alpha}) = |I_{\alpha}| = \lambda$.
\item[(b)] If $\alpha \neq \beta$ and $\varphi_{\alpha}: I_{\alpha} \to L$,
$\varphi_{\beta}: I_{\beta} \to L^{\prime}$ are invariant cofinal embeddings,
then $L \not \approx_{f} L^{\prime}$.
\end{enumerate}
\end{Lem}

\begin{proof}
Applying Solovay's Theorem, let $\{ S_{\tau} \mid \tau < \lambda \}$ 
be a partition of the stationary set
\[
S = \{ \delta < \lambda \mid \cf(\delta) = \omega_{1} \}
\]
into $\lambda$ pairwise disjoint stationary subsets. (For
example, see Jech \cite[7.6]{j}.) Fix some subset $X \subset \lambda$. 
Then for each $\alpha < \lambda$,
we define
\[
\lambda^{X}_{\alpha} =
\begin{cases}
\omega_{2},  &\text{if } \alpha \in \bigcup_{\tau \in X}S_{\tau} \\
\omega_{1},  &\text{otherwise;}
\end{cases}
\]
and we define the linear order
\[
I_{X} = \{ (\alpha, \beta) \mid \alpha < \lambda \text{ and }
\beta < \lambda^{X}_{\alpha} \}
\]
by setting $(\alpha_{1},\beta_{1}) < (\alpha_{2},\beta_{2})$ if and
only if either:
\begin{itemize}
\item $\alpha_{1} < \alpha_{2}$, or
\item $\alpha_{1} = \alpha_{2}$ and $\beta_{1} > \beta_{2}$.
\end{itemize} 
Suppose that $X \neq Y \subseteq \lambda$. Let $L$, $L^{\prime}$ be
linear orders and let $\varphi_{X}: I_{X} \to L$,
$\varphi_{Y}: I_{Y} \to L^{\prime}$ be invariant cofinal embeddings. 
Suppose that $L \approx_{f} L^{\prime}$ and let $\psi: M \to M^{\prime}$
be an isomorphism between the final segments $M$, $M^{\prime}$ of
$L$, $L^{\prime}$ respectively. For each $\delta < \lambda$, let
\[
M_{\delta} = \{ m \in M \mid m < \varphi_{X}(\gamma,0) 
\text{ for some } \gamma < \delta \}
\]
and 
\[
M^{\prime}_{\delta} = \{ m^{\prime} \in M^{\prime} \mid 
m^{\prime} < \varphi_{Y}(\gamma,0) \text{ for some } \gamma < \delta \}.
\]
Then there exists a club $C \subseteq \lambda$ such that
$\psi[M_{\delta}] = M^{\prime}_{\delta}$ for all $\delta \in C$. Without
loss of generality, we can suppose that there exists an ordinal
$\tau \in X \smallsetminus Y$. Choose $\delta \in C \cap S_{\tau}$
such that $M_{\delta} \neq \emptyset$. Since $\varphi_{X}$ is an
invariant embedding, it follows that 
$\ci(M \smallsetminus M_{\delta}) = \omega_{2}$. Similarly,
$\ci(M^{\prime} \smallsetminus M^{\prime}_{\delta}) = \omega_{1}$. 
But this is impossible since $\psi[M \smallsetminus M_{\delta}] =
M^{\prime} \smallsetminus M^{\prime}_{\delta}$.
\end{proof}

Note that in the statement of Theorem \ref{T:code}, $\kappa$ is not 
necessarily regular. In Section \ref{S:ultra}, we shall apply Theorem
\ref{T:code} in the case when $\kappa = 2^{\omega} > \omega_{1}$. 

\begin{Thm} \label{T:code}
If $\kappa > \omega_{1}$, then there exists a set 
$\{ J_{\alpha} \mid \alpha < 2^{\kappa} \}$ of linear orders satisfying the 
following conditions:
\begin{enumerate}
\item[(a)] $|J_{\alpha}| = \kappa$.
\item[(b)] $\ci(J_{\alpha}) = \cf(\kappa) + \omega_{2}$.
\item[(c)] If $\alpha \neq \beta$ and $\varphi_{\alpha}: J_{\alpha} \to L$,
$\varphi_{\beta}: J_{\beta} \to L^{\prime}$ are invariant coinitial embeddings,
then $L \not \approx_{i} L^{\prime}$.
\end{enumerate}
\end{Thm}

\begin{proof}
Let $\langle \kappa_{i} \mid i < \cf(\kappa) \rangle$ be a sequence
of regular cardinals with each $\kappa_{i} > \omega_{1}$ such that:
\begin{enumerate}
\item[(a)] if $\kappa$ is singular, then 
$\kappa = \sup_{i < \cf(\kappa)} \kappa_{i}$; and
\item[(b)] if $\kappa$ is regular, then $\kappa_{i} = \kappa$
for all $i < \cf(\kappa) = \kappa$.
\end{enumerate}
In either case, we have that 
$\prod_{i < \cf(\kappa)} 2^{\kappa_{i}} = 2^{\kappa}$. (For the
case when $\kappa$ is singular, see Jech \cite[6.5]{j}.)
Let $\theta = \cf(\kappa) + \omega_{2}$ and let 
$\{ S_{\tau} \mid \tau < \cf(\kappa) \}$ 
be a partition of the stationary set
\[
S = \{ \delta < \theta \mid \cf(\delta) = \omega_{1} \}
\]
into $\cf(\kappa)$ pairwise disjoint stationary subsets. 
Let $h: \theta \to \cf(\kappa)$ be the function defined by:
\begin{itemize}
\item $\delta \in S_{h(\delta)}$ for all $\delta \in S$; and
\item $h(\xi) = 0$ for all $\xi \in \theta \smallsetminus S$.
\end{itemize}
For each $i < \cf(\kappa)$, let 
$\{ I_{i, \alpha} \mid \alpha < 2^{\kappa_{i}} \}$ be the set
of linear orders of cardinality $\kappa_{i}$ given by 
Lemma \ref{L:code}. For any 
$\nu \in \prod_{i < \cf(\kappa)} 2^{\kappa_{i}}$, we define
the linear order
\[
J_{\nu} = \{ (\alpha, x ) \mid \alpha < \theta \text{ and }
x \in I_{h(\alpha), \nu(h(\alpha))} \}
\] 
by setting $(\alpha_{1}, x_{1}) < (\alpha_{2}, x_{2})$ if and 
only if either:
\begin{itemize}
\item $\alpha_{1} > \alpha_{2}$, or
\item $\alpha_{1} = \alpha_{2}$ and $x_{1} < x_{2}$.
\end{itemize} 
Arguing as in the proof of Lemma \ref{L:code}, we see that the set 
$\{ J_{\nu} \mid \nu \in \prod_{i < \cf(\kappa)} 2^{\kappa_{i}} \}$
of linear orders satisfies our requirements.
\end{proof}

Now suppose that $2^{\omega} = \kappa > \omega_{1}$ and
let $\{ J_{\alpha} \mid \alpha < 2^{\kappa} \}$ be the set of
linear orders given by Theorem \ref{T:code}. Then clearly
Theorem \ref{T:point} would follow if we could construct a set
$\{ \mathcal{D}_{\alpha} \mid \alpha < 2^{\kappa} \}$ of
nonprincipal ultrafilters over $\omega$ such that for each
$\alpha < 2^{\kappa}$, there exists an invariant coinitial
embedding
\[
\varphi_{\alpha} : J_{\alpha} \to 
(\omega^{\omega}/\mathcal{D}_{\alpha})^{*}.
\] 
Unfortunately this direct approach leads to serious technical 
difficulties which we have not yet been able to overcome. In
order to avoid these difficulties, in the next section, we
shall instead construct a set
$\{ \mathcal{D}_{\alpha} \mid \alpha < 2^{\kappa} \}$ of
nonprincipal ultrafilters over $\omega$ such that the
following condition is satisfied:
\begin{itemize}
\item For each $\alpha < 2^{\kappa}$ and each initial segment
$L$ of $(\omega^{\omega}/\mathcal{D}_{\alpha})^{*}$, there
exists an invariant embedding
$\varphi: \omega_{1} + J_{\alpha} \to L$.
\end{itemize}
(Here $\omega_{1} + J_{\alpha}$ is the linear order
consisting of a copy of the ordinal $\omega_{1}$ followed
by a copy of $J_{\alpha}$. In particular,
$( \omega_{1}, J_{\alpha})$ is an 
$(\omega_{1}, \cf(\kappa) + \omega_{2})$-cut of 
$\omega_{1} + J_{\alpha}$.)
Of course, this is not enough to ensure that
\[
( \omega^{\omega}/ \mathcal{D}_{\alpha} )^{*} \not \approx_{i}
( \omega^{\omega}/ \mathcal{D}_{\beta} )^{*} 
\]
for all $\beta \neq \alpha$, since the above condition
does not rule out the possibility that there also exists an 
invariant embedding
\[
\psi: \omega_{1} + J_{\beta} \to 
( \omega^{\omega}/ \mathcal{D}_{\alpha} )^{*} .
\]
Fix some $\alpha < 2^{\kappa}$ and let $C_{\alpha}$ be the set of 
$\beta < 2^{\kappa}$ such that there exists an invariant embedding
\[
\psi_{\beta}: \omega_{1} + J_{\beta} \to 
( \omega^{\omega}/ \mathcal{D}_{\alpha} )^{*} .
\]
For each $\beta \in C_{\alpha}$, 
let $(A_{\beta}, B_{\beta})$ be the
$(\omega_{1}, \cf(\kappa) + \omega_{2})$-cut of 
$( \omega^{\omega}/ \mathcal{D}_{\alpha} )^{*}$ defined by
\[
A_{\beta} = \{ g/\mathcal{D}_{\alpha} \in 
( \omega^{\omega}/ \mathcal{D}_{\alpha} )^{*}  \mid
g/\mathcal{D}_{\alpha} < \psi_{\beta}(t) 
\text{ for some } t \in \omega_{1} \}
\]
and
\[
B_{\beta} = \{ g/\mathcal{D}_{\alpha} \in 
( \omega^{\omega}/ \mathcal{D}_{\alpha} )^{*}  \mid
g/\mathcal{D}_{\alpha} > \psi_{\beta}(t) 
\text{ for some } t \in J_{\beta} \}.
\]
Then Theorem \ref{T:code} implies that 
$(A_{\beta}, B_{\beta}) \neq (A_{\gamma}, B_{\gamma})$ for
all $\beta \neq \gamma \in C_{\alpha}$. Since the following result 
implies that the number of $(\omega_{1}, \cf(\kappa) + \omega_{2})$-cuts 
of $( \omega^{\omega}/ \mathcal{D}_{\alpha} )^{*}$ is at
most $2^{\omega} = \kappa$, it follows that $|C_{\alpha}| \leq \kappa$.
This implies that there exists a subset
$W \subseteq 2^{\kappa}$ of cardinality $2^{\kappa}$ such that
\[
( \omega^{\omega}/ \mathcal{D}_{\alpha} )^{*} \not \approx_{i}
( \omega^{\omega}/ \mathcal{D}_{\beta} )^{*} 
\]
for all $\alpha \neq \beta \in W$. 

While the basic idea of Theorem \ref{T:cut} is implicitly
contained in Section VIII.0 of Shelah \cite{s} and 
Chapter III of Shelah \cite{s2}, the result does not
seem to have been explicitly stated anywhere in the literature.

\begin{Thm} \label{T:cut}
Suppose that $I$ is a linear order and that $\theta \neq \lambda$ are
regular cardinals. Then the number of $(\lambda, \theta)$-cuts of $I$ 
is at most $|I|$. 
\end{Thm}

\begin{proof}
We shall just consider the case when $\lambda < \theta$.
Suppose that $I$ is a counterexample of minimal cardinality and
let $\{ (A_{i},B_{i}) \mid i < |I|^{+} \}$ be a set of $|I|^{+}$
distinct $(\lambda, \theta)$-cuts of $I$. Let $\cf(|I|) = \kappa$ 
and express $I = \bigcup_{\gamma < \kappa}I_{\gamma}$ as  
a smooth strictly increasing union of substructures such that
$|I_{\gamma}| < |I|$ for all $\gamma < \kappa$.

First suppose that $\kappa \neq \theta$, $\lambda$. Then for
each $i < |I|^{+}$, there exists an ordinal $\gamma_{i} < \kappa$ 
such that $A_{i} \cap I_{\gamma_{i}}$ is cofinal in $A_{i}$ and
$B_{i} \cap I_{\gamma_{i}}$ is coinitial in $B_{i}$. It follows 
that there exists a subset $X \subseteq |I|^{+}$ of cardinality
$|I|^{+}$ and a fixed ordinal $\gamma < \kappa$ such that 
$\gamma_{i} = \gamma$ for all $i \in X$. But this means that
$\{ (A_{i} \cap I_{\gamma},B_{i} \cap I_{\gamma}) 
\mid i < X \}$ is a set of $|I|^{+}$ distinct 
$(\lambda, \theta)$-cuts of $I_{\gamma}$, which contradicts the
minimality of $|I|$.

Next suppose that $\kappa = \lambda$. Once again, for
each $i < |I|^{+}$, there exists an ordinal $\gamma_{i} < \kappa$ 
such that $B_{i} \cap I_{\gamma_{i}}$ is coinitial in $B_{i}$; and 
there exists a subset $X \subseteq |I|^{+}$ of cardinality
$|I|^{+}$ and a fixed ordinal $\gamma < \kappa$ such that 
$\gamma_{i} = \gamma$ for all $i \in X$. Arguing as in the previous 
paragraph, we can suppose that for each $i \in X$,
$A_{i} \cap I_{\gamma}$ is {\em not\/} cofinal in $A_{i}$.
For each $i \in X$, choose
an element $a_{i} \in A_{i} \smallsetminus I_{\gamma}$ such that
$s < a_{i} < t$ for all $s \in A_{i} \cap I_{\gamma}$ and
$t \in B_{i}$. Suppose that $i \neq j \in X$. Then we can suppose
that $B_{i} \varsubsetneq B_{j}$. Since $B_{j} \cap I_{\gamma}$ 
is coinitial in $B_{j}$, it follows that there exists an
element
\[
c \in (B_{j} \smallsetminus B_{i}) \cap I_{\gamma} \subseteq
A_{i} \cap I_{\gamma}.
\]
But this means that $a_{j} < c < a_{i}$ and so
$\{ a_{i} \mid i \in X \}$ is a set of $|I|^{+}$ distinct
elements of $I$, which is a
contradiction. A similar argument handles the case when
$\kappa = \theta$. This completes the proof of Theorem \ref{T:cut}.
\end{proof}

\section{Constructing ultrafilters} \label{S:ultra}

In this section, we shall prove Theorem \ref{T:point}. Our construction
of the required set $\{ \mathcal{D}_{\alpha} \mid \alpha < 2^{2^{\omega}} \}$ 
of nonprincipal ultrafilters makes use of the techniques developed in
Section VI.3 of Shelah \cite{s}. 

\begin{Def} \label{D:mod}
Let $D$ be a filter over $\omega$ and let
\[
I_{D} = \{ X \subseteq \omega \mid \omega \smallsetminus X \in D \}
\]
be the corresponding dual ideal. If $A$, $B \subseteq \omega$, then
we define
\[
A \subset B \mod D \quad \text{ if and only if } \quad
A \smallsetminus B \in I_{D}
\]
and
\[
A = B \mod D \quad \text{ if and only if } \quad
(A \smallsetminus B) \cup (B \smallsetminus A) \in I_{D}. 
\]
\end{Def}

\begin{Def} \label{D:ind1}
Suppose that $D$ is a filter over $\omega$ 
and that $\mathcal{G} \subseteq \omega^{\omega}$ is a family
of surjective functions. Then $\mathcal{G}$ is 
{\em independent mod $D$\/} if
for all distinct $g_{1}, \cdots, g_{\ell} \in \mathcal{G}$ and
all (not necessarily distinct) $j_{1}, \cdots , j_{\ell} \in \omega$,
\[
\{ n \in \omega \mid g_{k}(n) = j_{k} \text{ for all }
1 \leq k \leq \ell \} \neq \emptyset \mod D.
\]
(Of course, this condition implies that $D$ is a nontrivial
filter.)
\end{Def}

Suppose that $\mathcal{G}$ is independent mod $D$ and that
$|\mathcal{G}| = \kappa$. Let $I$ be any linear order of
cardinality $\kappa$ and suppose that 
$\mathcal{G} = \{ f_{t} \mid t \in I \}$ is indexed by the
elements of $I$. For each $s < t \in I$, let
\[
B_{s,t} = \{ n \in \omega \mid f_{s}(n) < f_{t}(n) \}.
\]  
Then it is easily checked that 
$D \cup \{ B_{s,t} \mid s < t \in I \}$ generates a nontrivial filter
$D^{+}$. (For example, see the proof of Lemma \ref{L:main}.)
It follows that if $\mathcal{D} \supseteq D^{+}$ is an
ultrafilter, then we can define an order-preserving map
$\varphi: I \to \omega^{\omega}/\mathcal{D}$ by
$\varphi(t) = f_{t}/\mathcal{D}$. However, if we wish 
$\varphi$ to be an invariant embedding, then we need to be able to
control the behaviour of arbitary elements 
$g/\mathcal{D} \in \omega^{\omega}/\mathcal{D}$. 
The next few paragraphs will introduce the techniques which
will enable us to accomplish this.

\begin{Def} \label{D:ind2}
Suppose that $\mathcal{G} \subseteq \omega^{\omega}$ is a family
of surjective functions. 
\begin{enumerate}
\item[(a)] $FI(\mathcal{G})$ is the set of functions $h$ satisfying
the following conditions:
\begin{enumerate}
\item[(i)] $\dom h$ is a finite subset of $\mathcal{G}$;
\item[(ii)] $\ran h \subset \omega$.
\end{enumerate}
\item[(b)] For each $h \in FI(\mathcal{G})$, let
\[
A_{h} = \{ n \in \omega \mid g(n) = h(g) \text{ for all }
g \in \dom h \}.
\]
\item[(c)] $FI_{s}(\mathcal{G}) = 
\{ A_{h} \mid h \in FI(\mathcal{G}) \}$.
\end{enumerate}
\end{Def}

\begin{Lem} \label{L:ind}
Suppose that $D$ is a filter over $\omega$ and that
$\mathcal{G} \subseteq \omega^{\omega}$ is a family
of surjective functions. 
\begin{enumerate}
\item[(a)] $\mathcal{G}$ is independent mod $D$ if and only if
$A_{h} \neq \emptyset$ mod $D$ for every $h \in FI(\mathcal{G})$. 
\item[(b)] If $\mathcal{G}$ is independent mod $D$, then there
exists a maximal filter $D^{*} \supseteq D$ modulo which
$\mathcal{G}$ is independent.
\item[(c)] If $\mathcal{G}$ is independent mod $D$ and
$X \subseteq \omega$, then there exists a finite subset
$\mathcal{F} \subseteq \mathcal{G}$ such that
$\mathcal{G} \smallsetminus \mathcal{F}$ is independent
modulo either the filter generated by $D \cup \{ X \}$
or the filter generated by $D \cup \{ \omega \smallsetminus X \}$.
\end{enumerate}
\end{Lem}

\begin{proof}
These are the statements of Claims 3.15(4), 3.15(3) and
3.3 from Shelah \cite[Chapter VI]{s}.
\end{proof}

Now suppose that $\mathcal{G} \subseteq \omega^{\omega}$ is a family
of surjective functions and that $D$ is a {\em maximal\/} filter over 
$\omega$ modulo which $\mathcal{G}$ is independent. Then
$\mathcal{A}$ is said to be a {\em partition mod $D$\/} if 
the following conditions are satisfied:
\begin{itemize}
\item $A \neq \emptyset \mod D$ for all $A \in \mathcal{A}$;
\item $A \cap A^{\prime} = \emptyset \mod D$ for all 
$A \neq A^{\prime} \in \mathcal{A}$;
\item for all $B \in \mathcal{P}(\omega)$ with $B \neq \emptyset \mod D$,
there exists $A \in \mathcal{A}$ such that $A \cap B \neq \emptyset \mod D$.
\end{itemize}
The subset $B \subseteq \omega$ is said to be {\em based on\/}
$\mathcal{A}$ if for every $A \in \mathcal{A}$, either
$A \subseteq B \mod D$ or $A \cap B = \emptyset \mod D$. 
By Claim 3.17(1) \cite[Chapter VI]{s}, for every subset $B \subseteq \omega$, 
there exists a partition $\mathcal{A}$ mod $D$ such that 
\begin{enumerate}
\item[(i)] $B$ is based on $\mathcal{A}$; and
\item[(ii)] $\mathcal{A} \subseteq FI_{s}(\mathcal{G})$. 
\end{enumerate}
Furthermore, by Claim 3.17(5) \cite[Chapter VI]{s}, $\mathcal{A}$ is
necessarily countable and so there exists a countable subset
$\mathcal{G}_{0} \subseteq \mathcal{G}$ such that
$\mathcal{A} \subseteq FI_{s}(\mathcal{G}_{0})$. In this case,
we say that $B$ is {\em supported by\/} $FI_{s}(\mathcal{G}_{0})$ mod $D$.

The next lemma summarises the properties of supports that we shall
require later in this section.

\begin{Lem} \label{L:supp}
Suppose that $\mathcal{G} \subseteq \omega^{\omega}$ is a family
of surjective functions and that $D$ is a maximal filter over $\omega$ 
modulo which $\mathcal{G}$ is independent.
\begin{enumerate}
\item[(a)] $FI_{s}(\mathcal{G})$ is dense mod $D$; i.e. for every
$B \subseteq \omega$ with $B \neq \emptyset \mod D$, there exists
$A_{h} \in FI_{s}(\mathcal{G})$ such that $A_{h} \subseteq B \mod D$. 
\item[(b)] For each $B \subseteq \omega$, there exists a countable
subset $\mathcal{G}_{0} \subseteq \mathcal{G}$ such that $B$ is
supported by $FI_{s}(\mathcal{G}_{0})$ mod $D$.
\item[(c)] Suppose that $\mathcal{G} = \mathcal{G}_{1} \sqcup \mathcal{G}_{2}$ 
and that $A \subseteq \omega$ is supported by $FI_{s}(\mathcal{G}_{1})$
mod $D$. If $h \in FI(\mathcal{G})$ and $A_{h} \subseteq A \mod D$, then 
$A_{h_{1}} \subseteq A \mod D$, where $h_{1} = h \res \mathcal{G}_{1}$.
\end{enumerate}
\end{Lem}

\begin{proof}
We have already discussed clause (b). Clauses (a) and (c) are the statements of 
Claims 3.17(1) and 3.17(4) from Shelah \cite[Chapter VI]{s}. 
\end{proof}

The following lemma will ensure that our construction
concentrates on the initial segments of 
$(\omega^{\omega}/\mathcal{D})^{*}$.

\begin{Lem} \label{L:coi}
Suppose that $\mathcal{G} \subseteq \omega^{\omega}$ is a family
of surjective functions and that $D$ is a maximal filter over $\omega$ 
modulo which $\mathcal{G}$ is independent. Suppose also that
$g \in \omega^{\omega}$ is a function such that
$\ell < g/D$ for every $\ell \in \omega$. Then
$f/D < g/D$ for every $f \in \mathcal{G}$.
\end{Lem}

\begin{proof}
This is Claim 3.19(1) from Shelah \cite[Chapter VI]{s}.
\end{proof}

Finally the next lemma is the key to our construction of the required set
of ultrafilters.

\begin{Lem} \label{L:main}
Suppose that $\mathcal{G} \sqcup \mathcal{G}^{*} \subseteq \omega^{\omega}$ 
is a family of surjective functions and that $D$ is a maximal filter over 
$\omega$ modulo which $\mathcal{G} \sqcup \mathcal{G}^{*}$ is independent. 
Suppose that $I$ is a linear order and that 
$\mathcal{G} = \{ f_{t} \mid t \in I \}$
is indexed by the elements of $I$. Then there exists a filter
$D^{+} \supseteq D$ over $\omega$ which satisfies the following
conditions:
\begin{enumerate}
\item[(a)] If $s \in I$ and $\ell \in \omega$, then $\ell < f_{s}/D^{+}$.
\item[(b)] If $s < t \in I$, then $f_{s}/D^{+} < f_{t}/D^{+}$.
\item[(c)] Suppose that $(I_{1}, I_{2})$ is a 
$(\lambda, \theta)$-cut of $I$ such that $\lambda$, $\theta > \omega$.
Then for every ultrafilter $\mathcal{U} \supseteq D^{+}$ over
$\omega$, there does {\em not\/} exist a function $g \in \omega^{\omega}$
such that 
\[
f_{s}/\mathcal{U} < g/\mathcal{U} < f_{t}/\mathcal{U}
\]
for all $s \in I_{1}$, $t \in I_{2}$.
\item[(d)] $D^{+}$ is a maximal filter over 
$\omega$ modulo which $\mathcal{G}^{*}$ is independent. 
\end{enumerate}
\end{Lem}

\begin{proof}
For each $t \in I$ and $\ell \in \omega$, let
\[
A_{\ell, t} = \{ n \in \omega \mid \ell < f_{t}(n) \}.
\] 
For each pair $s < t \in I$, let
\[
B_{s,t} = \{ n \in \omega \mid f_{s}(n) < f_{t}(n) \}. 
\]
Finally for each pair $r < s \in I$ and each function
$g \in \omega^{\omega}$ such that $g^{-1}(\ell)$ is supported by 
$FI_{s}( \mathcal{G}^{*} \sqcup \{ f_{t} \mid t \in I \smallsetminus [r,s] \})$
mod $D$ for all $\ell \in \omega$, let
\[
C_{g,r,s} = \{ n \in \omega \mid g(n) < f_{r}(n) \text{ or }
f_{s}(n) < g(n) \}.
\]
Let $E$ be the filter on $\omega$ generated by $D$, together with all of the sets
$A_{\ell, t}$, $B_{s,t}$, $C_{g,r,s}$ defined above. (Note that the next claim 
implies that $E$ is a nontrivial filter.)

\begin{Claim} \label{C:main}
If $h \in FI( \mathcal{G}^{*})$, then $A_{h} \neq \emptyset$ mod $E$.
\end{Claim}

Assuming Claim \ref{C:main}, we shall now complete the proof of 
Lemma \ref{L:main}. Applying Lemma \ref{L:ind}(a), $\mathcal{G}^{*}$ 
is independent mod $E$. Let $D^{+} \supseteq E$ be a maximal such filter.

Since $A_{\ell, t} \in D^{+}$ for each $t \in I$ and $\ell \in \omega$, 
it follows that $\ell < f_{s}/D^{+}$ and so clause \ref{L:main}(a) holds.
Similarly, since $B_{s,t} \in D^{+}$ for each $s < t \in I$, it follows
that clause \ref{L:main}(b) holds. Finally suppose that $(I_{1}, I_{2})$ 
is a $(\lambda, \theta)$-cut of $I$ such that $\lambda$, $\theta > \omega$
and that $g \in \omega^{\omega}$. By Lemma \ref{L:supp}(b), for each
$\ell \in \omega$, there exists a countable subset 
$\mathcal{G}_{\ell} \subseteq \mathcal{G} \sqcup \mathcal{G}^{*}$ 
such that $g^{-1}(\ell)$ is supported by $FI_{s}(\mathcal{G}_{\ell})$ 
mod $D$. Since $\lambda$, $\theta > \omega$, it follows that there
exist $r \in I_{1}$ and $s \in I_{2}$ such that $g^{-1}(\ell)$ is supported by 
$FI_{s}( \mathcal{G}^{*} \sqcup \{ f_{t} \mid t \in I \smallsetminus [r,s] \})$
mod $D$ for all $\ell \in \omega$ and hence $C_{g,r,s} \in D^{+}$. This
implies that if $\mathcal{U} \supseteq D^{+}$ is an ultrafilter over
$\omega$, then either $g/\mathcal{U} < f_{r}/\mathcal{U}$ or
$f_{s}/\mathcal{U} < g/\mathcal{U}$. Thus clause \ref{L:main}(c) also holds.

Thus it only remains to prove Claim \ref{C:main}. Suppose that
$h \in FI( \mathcal{G}^{*})$. Then it is enough to prove that
\[
A_{h} \cap \bigcap_{i \leq a}A_{\ell_{i},t_{i}} \cap
\bigcap_{i < j \leq a}B_{t_{i},t_{j}} \cap
\bigcap_{k \leq b}C_{g_{k}, r_{k},s_{k}} \neq \emptyset \mod D
\]
in the case when the following conditions are satisfied:
\begin{itemize}
\item $t_{0} < t_{1} < \cdots < t_{a}$.
\item If $k \leq b$, then $r_{k}$, $s_{k} \in \{ t_{i} \mid i \leq a \}$.
\end{itemize}
Let $\mathcal{T} = \{ f_{t_{i}} \mid i \leq a \}$.
We shall define a sequence of functions
$h_{m} \in FI(\mathcal{G} \sqcup \mathcal{G}^{*})$ inductively
for $m \in \omega$ so
that the following conditions are satisfied:
\begin{enumerate}
\item[(1)] $h_{0} = h$ and $h_{m} \subseteq h_{m+1}$.
\item[(2)] $\dom h_{m} \cap \mathcal{T} = \emptyset$.
\item[(3)] If $h^{*} \in FI(\mathcal{T})$ and $k \leq b$, 
then one of the following occurs for almost all $m$,
\begin{enumerate}
\item[(i)] there exists $\ell \in \omega$ such that
$A_{h_{m} \cup h^{*}} \subseteq g^{-1}_{k}(\ell) \mod D$; or
\item[(ii)] $A_{h_{m} \cup h^{*}} \cap g^{-1}_{k}(\ell) =
\emptyset \mod D$ for all $\ell \in \omega$.
\end{enumerate}
\end{enumerate}
(Clearly if (i) occurs, then there exists a {\em fixed\/} $\ell$ such
that $A_{h_{m} \cup h^{*}} \subseteq g^{-1}_{k}(\ell) \mod D$
for almost all $m$.) To see that the induction can be carried out, first fix
an enumeration of the countably many pairs $h^{*}$, $k$ that
must be dealt with. Now suppose that $h_{m}$ has been 
defined and that we must next deal with the pair
$h^{*}$, $k$. There are two cases to consider. First suppose
that there exists $\ell \in \omega$ such that
$A_{h_{m} \cup h^{*}} \cap g^{-1}_{k}(\ell) \neq \emptyset \mod D$.
By Lemma \ref{L:supp}(a), there exists 
$\tilde{h} \in FI(\mathcal{G} \sqcup \mathcal{G}^{*})$ such that
\[
A_{\tilde{h}} \subseteq A_{h_{m} \cup h^{*}} \cap g^{-1}_{k}(\ell) 
\mod D.
\] 
Clearly we must have that $h_{m} \cup h^{*} \subseteq \tilde{h}$;
and in this case, we set 
\[
h_{m+1} = \tilde{h} \res (( \mathcal{G} \sqcup \mathcal{G}^{*} )
\smallsetminus \mathcal{T}).
\]
Otherwise, we must have that $A_{h_{m} \cup h^{*}} \cap g^{-1}_{k}(\ell) =
\emptyset \mod D$ for all $\ell \in \omega$; and in this case,
we set $h_{m+1} = h_{m}$. 

Now fix some $k \leq b$. Let $r_{k} = t_{\iota(k)}$ and 
$s_{k} = t_{\tau(k)}$. Let 
\[
\mathcal{T}_{k} = \{ f_{t_{i}} \mid i \notin [ \iota(k), \tau(k)] \}.
\] 
Suppose that $h^{*} \in FI(\mathcal{T})$. Then for all
sufficiently large $m$, either:
\begin{enumerate}
\item[(i)] there exists $\ell \in \omega$ such that
$A_{h_{m} \cup h^{*}} \subseteq g^{-1}_{k}(\ell) \mod D$; or
\item[(ii)] $A_{h_{m} \cup h^{*}} \cap g^{-1}_{k}(\ell) =
\emptyset \mod D$ for all $\ell \in \omega$. 
\end{enumerate}
First suppose that (i) holds. 
Since $g^{-1}_{k}(\ell)$ is supported by 
\[
FI_{s}( \mathcal{G}^{*} \sqcup \{ f_{t} \mid t \in I \smallsetminus
[r_{k},s_{k}] \}),
\]
Lemma \ref{L:supp}(c) implies that: 
\begin{enumerate}
\item[$(i)^{\prime}$] there exists $\ell \in \omega$ such that 
$A_{h_{m} \cup (h^{*} \res \mathcal{T}_{k})} \subseteq g^{-1}_{k}(\ell)$ 
mod $D$ for almost all $m$.
\end{enumerate}
Similarly, if (ii) holds, then Lemma \ref{L:supp}(c) implies that:
\begin{enumerate}
\item[$(ii)^{\prime}$] for almost all $m$,
$A_{h_{m} \cup (h^{*} \res \mathcal{T}_{k})} \cap g^{-1}_{k}(\ell) = \emptyset$
mod $D$ for all $\ell \in \omega$. 
\end{enumerate}
Let $\psi_{k}: \omega^{\mathcal{T}_{k}} \to \omega \cup \{ \infty \}$ 
be the function defined by
\[
\psi_{k}( h^{*} \res \mathcal{T}_{k} ) =
\begin{cases}
\ell     &\text{if $(i)^{\prime}$ holds;}  \\
\infty   &\text{if $(ii)^{\prime}$ holds.}  
\end{cases}
\]
For each infinite set $W \subseteq \omega$, let $W^{(\mathcal{T})}$ be the set of  
functions $h^{*}: \mathcal{T} \to W$ such that
$h^{*}(f_{t_{i}}) < h^{*}(f_{t_{j}})$ for all $i < j \leq a$; 
and for each $k \leq b$, let
$\varphi_{k}: \omega^{(\mathcal{T})} \to 3$ be the function defined by
\[
\varphi_{k}(h^{*}) =
\begin{cases}
0     &\text{if $\psi_{k}(h^{*} \res \mathcal{T}_{k}) < h^{*}(f_{t_{\iota(k)}})$;} \\
1     &\text{if $\psi_{k}(h^{*} \res \mathcal{T}_{k}) > h^{*}(f_{t_{\tau(k)}})$;} \\
2     &\text{otherwise.}
\end{cases}
\]
By Ramsey's Theorem, there exists an infinite set $W \subseteq \omega$
such that $\varphi_{k} \res W^{(\mathcal{T})}$ is a constant function for
all $k \leq b$. 

\begin{Claim} \label{C:ramsey}
For each $k \leq b$, either $\varphi_{k} \res W^{(\mathcal{T})} \equiv 0$ or
$\varphi_{k} \res W^{(\mathcal{T})} \equiv 1$.  
\end{Claim}

\begin{proof}
Suppose that $\varphi_{k} \res W^{(\mathcal{T})} \equiv 2$, so that
\[
h^{*}(f_{t_{\iota(k)}}) \leq \psi_{k}(h^{*} \res \mathcal{T}_{k})
\leq h^{*}(f_{t_{\tau(k)}})
\]
for all $h^{*} \in W^{(\mathcal{T})}$. Let
$|\{ j \mid \iota(k) \leq j \leq \tau(k) \}| = p$ and let
$h^{\prime}: \mathcal{T}_{k} \to W$ be a strictly increasing function 
such that 
\[
|\{ w \in W \mid h^{\prime}(f_{t_{\iota(k) - 1}}) < w <
h^{\prime}(f_{t_{\tau(k) +1}}) \}| = 2p.
\]
Then we can extend $h^{\prime}$ to a function
$h^{*} \in W^{(\mathcal{T})}$ such that either
\[
\psi_{k}(h^{*} \res \mathcal{T}_{k}) = \psi_{k}(h^{\prime}) < 
h^{*}(f_{t_{\iota(k)}})
\]
or 
\[
\psi_{k}(h^{*} \res \mathcal{T}_{k}) = \psi_{k}(h^{\prime}) > 
h^{*}(f_{t_{\tau(k)}}),
\]
which is a contradiction.
\end{proof}

Choose an increasing sequence
$j_{0} < j_{1} < \cdots < j_{a}$ of elements of $W$ such that
$j_{i} > \ell_{i}$ for each $i \leq a$; and let
$h^{*} \in FI(\mathcal{T})$ be the function defined by
$h^{*}(f_{t_{i}}) = j_{i}$ for each $i \leq a$. To complete the proof of
Claim \ref{C:main}, it is enough to show that for almost
all $m$,
\[
A_{h_{m} \cup h^{*}} \subseteq 
A_{h} \cap \bigcap_{i \leq a}A_{\ell_{i},t_{i}} \cap
\bigcap_{i < j \leq a}B_{t_{i},t_{j}} \cap
\bigcap_{k \leq b}C_{g_{k}, r_{k},s_{k}} \mod D.
\]
It is clear that $A_{h_{m}} \subseteq A_{h}$ for all $m$ and that
\[
A_{h^{*}} \subseteq \bigcap_{i \leq a}A_{\ell_{i},t_{i}} \cap
\bigcap_{i < j \leq a}B_{t_{i},t_{j}}. 
\]
Finally let $k \leq b$. If $\varphi_{k}(h^{*}) = 0$, 
then $\psi_{k}(h^{*} \res \mathcal{T}_{k}) < h^{*}(f_{t_{\iota(k)}})
< \infty$ and so for almost all $m$,
\[
A_{h_{m} \cup h^{*}} \subseteq \{ n \in \omega \mid
g_{k}(n) = \psi_{k}(h^{*} \res \mathcal{T}_{k}) < 
h^{*}(f_{t_{\iota(k)}}) = f_{t_{\iota(k)}}(n) \} \mod D.
\]
Similarly, if $\varphi_{k}(h^{*}) = 1$ and
$\psi_{k}(h^{*} \res \mathcal{T}_{k}) < \infty$,
then for almost all $m$,
\[
A_{h_{m} \cup h^{*}} \subseteq \{ n \in \omega \mid
f_{t_{\tau(k)}}(n) = h^{*}(f_{t_{\tau(k)}}) <
\psi_{k}(h^{*} \res \mathcal{T}_{k}) = g_{k}(n) \} \mod D.
\]
On the other hand, if $\varphi_{k}(h^{*}) = 1$ and
$\psi_{k}(h^{*} \res \mathcal{T}_{k}) = \infty$, then 
for almost all $m$,
\[
A_{h_{m} \cup h^{*}} \cap 
\bigcup_{\ell \leq h^{*}(f_{t_{\tau(k)}})}g^{-1}_{k}(\ell)
= \emptyset \mod D
\]
and so
\[
A_{h_{m} \cup h^{*}} \subseteq \{ n \in \omega \mid
f_{t_{\tau(k)}}(n) = h^{*}(f_{t_{\tau(k)}}) <
g_{k}(n) \} \mod D.
\]
Hence, in every case, we have that for almost all $m$,
\[
A_{h_{m} \cup h^{*}} \subseteq C_{g_{k}, r_{k},s_{k}} \mod D.
\]
This completes the proof of Lemma \ref{L:main}.
\end{proof}

It is now straightforward to construct a set of ultrafilters
satisfying the conclusion of Theorem \ref{T:point}.

\begin{proof}[Proof of Theorem \ref{T:point}]
Suppose that $2^{\omega} = \kappa > \omega_{1}$.
Let $\{ J_{\alpha} \mid \alpha < 2^{\kappa} \}$ be the set of
linear orders given by Theorem \ref{T:code}; and for each $\alpha < \kappa$,
let $I_{\alpha} = \omega_{1} + J_{\alpha}$. Fix some $\alpha < \kappa$.
To simplify notation, let $I_{\alpha} = I$. Then the corresponding
ultrafilter $\mathcal{D}_{\alpha} = \mathcal{D}$ is constructed as
follows.

Let $F_{0} = \{ X \subseteq \omega \mid |\omega \smallsetminus X| < \omega \}$ 
be the Fr\'{e}chet filter over $\omega$. By Theorem 1.5(1) of
Shelah \cite[Appendix]{s}, there exists a family 
$\mathcal{G} \subseteq \omega^{\omega}$ of surjective functions
of cardinality $\kappa = 2^{\omega}$
such that $\mathcal{G}$ is independent mod $F_{0}$. Let 
$\mathcal{P}(\omega) = \{ X_{\mu} \mid \mu < \kappa \}$ be an enumeration of the
powerset of $\omega$ and ``enumerate'' $\mathcal{G}$ as
$\{ f_{\xi}^{\mu} \mid \mu, \xi < \kappa \}$. We shall define by induction
on $\mu < \kappa$
\begin{itemize}
\item a decreasing sequence of subsets $\mathcal{G}_{\mu} \subseteq 
\{ f_{\xi}^{\nu} \mid \xi < \kappa \text { and }
\mu \leq \nu < \kappa \}$; and
\item an increasing sequence of filters $D_{\mu}$ over $\omega$
\end{itemize}
such that the following conditions are satisfied:
\begin{enumerate}
\item[(a)] $\mathcal{G}_{0} = \mathcal{G}$.
\item[(b)] $|\{ f_{\xi}^{\nu} \mid \xi < \kappa \text { and }
\mu \leq \nu < \kappa \} \smallsetminus \mathcal{G}_{\mu}|
\leq |\mu| + \omega$.
\item[(c)] $D_{\mu}$ is a maximal filter modulo which $\mathcal{G}_{\mu}$
is independent.
\item[(d)] Either $X_{\mu} \in D_{\mu + 1}$ or 
$\omega \smallsetminus X_{\mu} \in D_{\mu + 1}$.
\end{enumerate}
When $\mu = 0$, we let $D_{0} \supseteq F_{0}$ be a maximal filter
modulo which $\mathcal{G}_{0} = \mathcal{G}$ is independent. If $\mu$
is a limit ordinal, then we define
$\mathcal{G}_{\mu} = \bigcap_{\nu < \mu}\mathcal{G}_{\nu}$ and let
$D_{\mu} \supseteq \bigcup_{\nu < \mu}D_{\nu}$ be a maximal filter
modulo which $\mathcal{G}_{\mu}$ is independent. Finally suppose that
$\mu = \nu +1$. By Lemma \ref{L:ind}, there exists a finite subset
$\mathcal{F}_{\nu} \subseteq \mathcal{G}_{\nu}$ such that
$\mathcal{G}_{\nu} \smallsetminus \mathcal{F}_{\nu}$ is independent
modulo either the filter generated by $D_{\nu} \cup \{ X_{\nu} \}$
or the filter generated by $D_{\nu} \cup \{ \omega \smallsetminus X_{\nu} \}$.
Without loss of generality, suppose that 
$\mathcal{G}_{\nu} \smallsetminus \mathcal{F}_{\nu}$ is independent
modulo the filter $D_{\nu}^{\prime}$ generated by 
$D_{\nu} \cup \{ X_{\nu} \}$; and let 
$E_{\nu} \supseteq D_{\nu}^{\prime}$ be a maximal filter modulo
which $\mathcal{G}_{\nu} \smallsetminus \mathcal{F}_{\nu}$ is independent.
Let 
\[
\mathcal{G}_{\mu} = 
\{ f^{\tau}_{\xi} \in \mathcal{G}_{\nu} \smallsetminus \mathcal{F}_{\nu}
\mid \mu \leq \tau < \kappa \}.
\]
Note that 
\[
\mathcal{H}_{\nu} = 
\{ f^{\tau}_{\xi} \in \mathcal{G}_{\nu} \smallsetminus \mathcal{F}_{\nu}
\mid \tau = \nu \}
\]
has cardinality $\kappa$. Hence we can re-index $\mathcal{H}_{\nu}$
as $\mathcal{H}_{\nu} = \{ f^{\nu}_{t} \mid t \in I \}$. 
By Lemma \ref{L:main}, there exists a filter $D_{\mu} \supseteq E_{\nu}$
which satisfies the following conditions:
\begin{enumerate}
\item[(1)] If $s \in I$ and $\ell \in \omega$, then $\ell < f^{\nu}_{s}/D_{\mu}$.
\item[(2)] If $s < t \in I$, then $f^{\nu}_{s}/D_{\mu} < f^{\nu}_{t}/D_{\mu}$.
\item[(3)] Suppose that $(I_{1}, I_{2})$ is a 
$(\lambda, \theta)$-cut of $I$ such that $\lambda$, $\theta > \omega$.
Then for every ultrafilter $\mathcal{U} \supseteq D_{\mu}$ over
$\omega$, there does {\em not\/} exist a function $g \in \omega^{\omega}$
such that 
\[
f^{\nu}_{s}/\mathcal{U} < g/\mathcal{U} < f^{\nu}_{t}/\mathcal{U}
\]
for all $s \in I_{1}$, $t \in I_{2}$.
\item[(4)] $D_{\mu}$ is a maximal filter over 
$\omega$ modulo which $\mathcal{G}_{\mu}$ is independent. 
\end{enumerate}
Finally let $\mathcal{D} = \bigcup_{\mu < \kappa}D_{\mu}$. By clause (d),
$\mathcal{D}$ is an ultrafilter.

\begin{Claim} \label{C:coi}
If $L$ is an initial segment of
$(\omega^{\omega}/\mathcal{D})^{*}$, then there exists
$\mu < \kappa$ such that
$\{ f^{\mu}_{t}/\mathcal{D} \mid t \in I \} \subseteq L$.
\end{Claim}

\begin{proof}[Proof of Claim \ref{C:coi}]
Let $g/\mathcal{D} \in L$. Then there exists $\mu < \kappa$ such that
\[
A_{\ell} = \{ n \in \omega \mid \ell < g(n) \} \in D_{\mu}
\]
for all $\ell \in \omega$. Since $D_{\mu}$ is a maximal filter modulo 
which $\mathcal{G}_{\mu}$ is independent, Lemma \ref{L:coi} implies
that $f/D_{\mu} < g/D_{\mu}$ for all $f \in \mathcal{G}_{\mu}$.
Hence $\{ f^{\mu}_{t}/\mathcal{D} \mid t \in I \} \subseteq L$.
\end{proof}

From now on, it is necessary to write $\mathcal{D}_{\alpha}$, 
$I_{\alpha}$, etc. 

\begin{Claim} \label{C:many}
Fix some $\alpha < 2^{\kappa}$. Then the set
\[
E_{\alpha} = \{ \beta < 2^{\kappa} \mid ( \omega^{\omega}/ \mathcal{D}_{\alpha} )^{*} 
\approx_{i} ( \omega^{\omega}/ \mathcal{D}_{\beta} )^{*} \}
\]
has cardinality at most $\kappa$.
\end{Claim}

\begin{proof}[Proof of Claim \ref{C:many}]
Suppose that $|E_{\alpha}| \geq \kappa^{+}$. For each $\beta \in E_{\alpha}$,
let $L_{\beta}$, $M_{\beta}$ be initial segments of
$( \omega^{\omega}/ \mathcal{D}_{\beta} )^{*}$,
$( \omega^{\omega}/ \mathcal{D}_{\alpha} )^{*}$ respectively such that
there exists an isomorphism $\varphi_{\beta}: L_{\beta} \to M_{\beta}$.
By Claim \ref{C:coi}, for each $\beta \in E_{\alpha}$, there exists
$\mu_{\beta} < \kappa$ such that
\[
R_{\beta} = \{ f^{\mu_{\beta}}_{t}/\mathcal{D}_{\beta} \mid t \in I_{\beta} \} 
\subseteq L_{\beta}.
\]
Recall that $I_{\beta} = \omega_{1} + J_{\beta}$. Let
$S_{\beta} = \{ f^{\mu_{\beta}}_{t}/\mathcal{D}_{\beta} \mid t \in \omega_{1} \}$ 
and
$T_{\beta} = \{ f^{\mu_{\beta}}_{t}/\mathcal{D}_{\beta} \mid t \in J_{\beta} \}$. 
Let $\theta = \cf(\kappa) + \omega_{2}$. Then
$( \varphi_{\beta}[S_{\beta}], \varphi_{\beta}[T_{\beta}] )$ determines the
$(\omega_{1}, \theta )$-cut $(A_{\beta}, B_{\beta})$ of
$( \omega^{\omega}/ \mathcal{D}_{\alpha} )^{*}$ defined by
\[
A_{\beta} = \{ g/\mathcal{D}_{\alpha} \in ( \omega^{\omega}/ \mathcal{D}_{\alpha} )^{*}
\mid g/\mathcal{D}_{\alpha} < \varphi_{\beta}(s) \text{ for some }
s \in S_{\beta} \}
\]
and 
\[
B_{\beta} = \{ g/\mathcal{D}_{\alpha} \in ( \omega^{\omega}/ \mathcal{D}_{\alpha} )^{*}
\mid g/\mathcal{D}_{\alpha} > \varphi_{\beta}(t) \text{ for some }
t \in T_{\beta} \}.
\]
By Theorem \ref{T:cut}, there exist $\beta \neq \gamma \in E_{\alpha}$
such that $(A_{\beta}, B_{\beta}) = (A_{\gamma}, B_{\gamma})$. But this
is impossible, since we can define invariant coinitial embeddings of $J_{\beta}$,
$J_{\gamma}$ into $B_{\beta} = B_{\gamma}$ by
$c \mapsto \varphi_{\beta}(f^{\mu_{\beta}}_{c})$ and
$d \mapsto \varphi_{\gamma}(f^{\mu_{\gamma}}_{d})$ respectively, which
contradicts Theorem \ref{T:code}.
\end{proof}
Clearly Claim \ref{C:many} implies that there exists a subset
$W \subseteq 2^{\kappa}$ of cardinality $2^{\kappa}$ such that
\[
( \omega^{\omega}/ \mathcal{D}_{\alpha} )^{*} \not \approx_{i}
( \omega^{\omega}/ \mathcal{D}_{\beta} )^{*} 
\]
for all $\alpha \neq \beta \in W$. This completes the proof
of Theorem \ref{T:point}.
\end{proof}

\section{Asymptotic cones under $CH$} \label{S:ch}

In this section, we shall prove Theorem \ref{T:ch}. 
Let $\Gamma$ be an infinite finitely generated 
group and let $d$ be the word metric with respect to some finite 
generating set. Then the main point is that each asymptotic cone 
$\Con_{\mathcal{D}}(\Gamma)$ can be uniformly constructed from an 
associated ultraproduct
\[
\prod \mathcal{M}_{n} / \mathcal{D},
\]
where each $\mathcal{M}_{n}$ is a suitable countable structure
for a fixed countable first-order language $\mathcal{L}$. If
$CH$ holds, then $\prod \mathcal{M}_{n} / \mathcal{D}$ is a saturated
structure of cardinality $\omega_{1}$ and hence is determined up to
isomorphism by its complete first-order theory $T_{\mathcal{D}}$. 
Consequently, if $CH$ holds, then
since there are at most $2^{\omega}$ 
possibilities for $T_{\mathcal{D}}$, there are also at most $2^{\omega}$ 
possibilities for $\prod \mathcal{M}_{n} / \mathcal{D}$ 
and hence also for $\Con_{\mathcal{D}}(\Gamma)$.

\begin{Def} \label{D:language}
Let $\mathcal{L}$ be the first-order language consisting of the following
symbols:
\begin{enumerate}
\item[(a)] the binary relation symbol $R_{q}$ for each 
$0 < q \in \mathbb{Q}$\:; and
\item[(b)] the constant symbol $e$.
\end{enumerate}
\end{Def}

\begin{Def} \label{D:structure}
For each $n \geq 1$, $\mathcal{M}_{n}$ is the $\mathcal{L}$-structure
with universe $\Gamma$ such that:
\begin{enumerate}
\item[(a)] $\mathcal{M}_{n} \vDash R_{q}(x,y)$ if and only if
$d(x,y) \leq qn$; and
\item[(b)] $e$ is the identity element of $\Gamma$.
\end{enumerate}
\end{Def}

\begin{Def} \label{D:foi}
For each nonprincipal ultrafilter $\mathcal{D}$ over $\omega$, 
let $T_{\mathcal{D}}$ be the complete first-order theory of
$\prod \mathcal{M}_{n} / \mathcal{D}$.
\end{Def}

\begin{Thm} \label{T:foi}
If $\mathcal{D}$, $\mathcal{D}^{\prime}$ are nonprincipal ultrafilters
over $\omega$, then the following are equivalent.
\begin{enumerate}
\item[(a)] $\prod \mathcal{M}_{n} / \mathcal{D} \cong
\prod \mathcal{M}_{n} / \mathcal{D}^{\prime}$.
\item[(b)] $\Con_{\mathcal{D}}(\Gamma)$ is isometric to
$\Con_{\mathcal{D}^{\prime}}(\Gamma)$. 
\end{enumerate}
\end{Thm}

\begin{proof}
We shall begin by describing how the asymptotic cone $\Con_{\mathcal{D}}(\Gamma)$
can be uniformly constructed from the ultraproduct
\[
\prod \mathcal{M}_{n} / \mathcal{D} = \langle X ; R_{q}, e \rangle.
\]
First define $\mathcal{M}^{0}_{\mathcal{D}}$ to be the set of those $x \in X$ such that
\[
\prod \mathcal{M}_{n} / \mathcal{D}  \vDash R_{q}(x,e)
\]
for some $q > 0$. Next define an equivalence relation $\approx$ on
$\mathcal{M}^{0}_{\mathcal{D}}$ by
\[
x \approx y \quad \text{ if and only if } \quad 
\prod \mathcal{M}_{n} / \mathcal{D}  \vDash R_{q}(x,y) \text{ for all } q > 0.
\]
For each $x \in \mathcal{M}^{0}_{\mathcal{D}}$, let $\langle x \rangle$ denote 
the corresponding $\approx$-class and let 
\[
C_{\mathcal{D}} = \{ \langle x \rangle \mid 
x \in \mathcal{M}^{0}_{\mathcal{D}} \}.
\]
Then we can define a metric $d_{\mathcal{D}}$ on $C_{\mathcal{D}}$ by
\[
d_{\mathcal{D}}(\langle x \rangle,\langle y \rangle) = 
\inf \{ q \mid \prod \mathcal{M}_{n} / \mathcal{D}  \vDash R_{q}(x,y) \}.
\]
It is easily checked that $\langle C_{\mathcal{D}}, d_{\mathcal{D}} \rangle$
is isometric to the asymptotic cone $\Con_{\mathcal{D}}(\Gamma)$. Consequently,
if $\prod \mathcal{M}_{n} / \mathcal{D} \cong
\prod \mathcal{M}_{n} / \mathcal{D}^{\prime}$, then
$\Con_{\mathcal{D}}(\Gamma)$ is isometric to
$\Con_{\mathcal{D}^{\prime}}(\Gamma)$. 

It is also easily checked that $\prod \mathcal{M}_{n} / \mathcal{D}$ consists
of $2^{\omega}$ disjoint isomorphic copies of $\mathcal{M}^{0}_{\mathcal{D}}$ 
and that each $\approx$-class has cardinality $2^{\omega}$. It follows that
any isometry between $\Con_{\mathcal{D}}(\Gamma)$ and 
$\Con_{\mathcal{D}^{\prime}}(\Gamma)$ can be lifted to a corresponding
isomorphism between $\prod \mathcal{M}_{n} / \mathcal{D}$ and
$\prod \mathcal{M}_{n} / \mathcal{D}^{\prime}$. 
\end{proof}

\begin{Cor} \label{C:foi}
Assume $CH$. If $\mathcal{D}$, $\mathcal{D}^{\prime}$ are nonprincipal 
ultrafilters over $\omega$, then the following are equivalent.
\begin{enumerate}
\item[(a)] $T_{\mathcal{D}} = T_{\mathcal{D}^{\prime}}$.
\item[(b)] $\Con_{\mathcal{D}}(\Gamma)$ is isometric to
$\Con_{\mathcal{D}^{\prime}}(\Gamma)$. 
\end{enumerate}
\end{Cor}

\begin{proof}
As we mentioned earlier, if
$CH$ holds, then $\prod \mathcal{M}_{n} / \mathcal{D}$ is a saturated
structure of cardinality $\omega_{1}$ and hence is determined up to
isomorphism by its complete first-order theory $T_{\mathcal{D}}$. 
(For example, see Chang-Keisler \cite{ck}.) 
\end{proof}

This completes the proof of Theorem \ref{T:ch}.

\end{document}